\documentclass[preprint,sort&compress,12pt]{elsarticle}
\usepackage{mathrsfs}
\usepackage{mathrsfs}
\usepackage{amsmath}
\usepackage{mathrsfs}
\usepackage{stmaryrd}
\usepackage{mathrsfs}
\usepackage{bm}
\usepackage{amsmath}
\usepackage{amsfonts}
\usepackage{amssymb}
\usepackage{amsthm}
\usepackage[RGB]{xcolor}
\journal{\quad}
\newcommand{\mf}{\mathbf}
\newcommand{\mm}{\mathrm}

\begin{document}

\begin{frontmatter}
\title{
On Multi-dimensional Compressible Flows of
Nematic Liquid \\ Crystals with Large Initial Energy in a Bounded Domain}
% \tnotetext[S]{The research of Fei Jiang was supported by the NSFC (Grant No. 11101044),
% the research of Song Jiang by NSFC (Grant No. 40890154) and the National Basic
% Research Program under the Grant 2011CB309705, and the research of
% Dehua Wang by the National Science Foundation under Grant DMS-0906160 and the Office of Naval
% Research under Grant N00014-07-1-0668.}
%%%%%
\author[FJ]{Fei Jiang}  %\corref{cor1}}
\ead{jiangfei0591@163.com}
\author[SJ]{Song Jiang}
\ead{jiang@iapcm.ac.cn}
% \cortext[cor1]{Corresponding
%author: Tel +86 15001201710.}
\author[DH]{Dehua Wang}
\ead{dwang@math.pitt.edu}
\address[FJ]{College of Mathematics and Computer Science, Fuzhou University, Fuzhou, 350108, China.}
\address[SJ]{Institute of Applied Physics and Computational Mathematics, Beijing, 100088, China.}
\address[DH]{Department of Mathematics, University of Pittsburgh, Pittsburgh, PA, 15260, USA.}

\begin{abstract}
We study the global existence of weak solutions to a multi-dimensional simplified
Ericksen-Leslie system for compressible flows of nematic liquid crystals with large initial energy
in a bounded domain $\Omega\subset \mathbb{R}^N$, where $N=2$ or $3$.
By exploiting a maximum principle, Nirenberg's interpolation inequality and a smallness condition
imposed on the $N$-th component of initial direction field $\mf{d}_0$ to overcome the difficulties
induced by the supercritical nonlinearity $|\nabla{\mathbf d}|^2{\mathbf d}$ in the
equations of angular momentum, and then adapting a modified three-dimensional
 approximation scheme and the weak convergence arguments for the compressible
Navier-Stokes equations, we establish the global existence of weak solutions to the initial-boundary
problem with large initial energy and without any smallness condition on the initial density and velocity.
\iffalse
We
\fi
\end{abstract}

\begin{keyword} Liquid crystals, compressible flows, weak solutions, weak convergence arguments.\\[0.5em]
\MSC[2000] 35Q35\sep 76D03.
%(2000 is the default)
\end{keyword}
\end{frontmatter}

% \linenumbers

\newtheorem{thm}{Theorem}[section]
\newtheorem{lem}{Lemma}[section]
\newtheorem{pro}{Proposition}[section]
\newtheorem{cor}{Corollary}[section]
\newproof{pf}{Proof}
\newdefinition{rem}{Remark}[section]
\newtheorem{definition}{Definition}[section]
\newcommand{\red}{\color{red}}

\section{Introduction}\label{sec:01}
\label{Intro} \numberwithin{equation}{section}

We study the global existence of weak solutions to the following multi-dimensional simplified version
of the Ericksen-Leslie model in a  bounded domain $\Omega\subset{\mathbb R}^N$ which describes
the motion of a compressible flow of nematic liquid crystals:
\begin{eqnarray}
&& \label{0101}\partial_t\rho+\mathrm{div}(\rho\mathbf{v})=0,\\
 &&  \label{0102} \partial_t
(\rho\mathbf{v})+\mm{div}(\rho\mathbf{v}\otimes \mathbf{v})+ \nabla
P=\mu \Delta\mathbf{v}+(\mu+\lambda)\nabla
\mm{div}\mf{v}\nonumber\\
&&\qquad \qquad \qquad \qquad  \qquad \qquad  \quad \qquad
-\nu\mm{div}\Big( \nabla \mathbf{d}\odot\nabla
\mathbf{d}-\frac{1}{2}|\nabla \mf{d}|^2\mathbb{I}\Big),  \\
&&  \label{0103}\partial_t\mathbf{d}+\mathbf{v}\cdot \nabla
\mathbf{d}= \theta(\Delta \mathbf{d}+|\nabla \mathbf{d}|^2\mf{d}),
\end{eqnarray} with initial conditions:
\begin{equation} \label{0105}
\rho(\mathbf{x},0)=\rho_0(\mathbf{x}),\ \
\mathbf{d}(\mathbf{x},0)=\mathbf{d}_0(\mathbf{x}),\ \ (\rho
\mathbf{v})(\mathbf{x},0)=\mathbf{m}_0(\mathbf{x})\quad\mbox{ in }\Omega ,
\end{equation}
and boundary conditions:
\begin{equation} \label{0104}
       \mf{n}\cdot\nabla{\mathbf{d}(\mathbf{x},t)}=\mf{0}, \
\mathbf{v}(\mathbf{x},t)=\mathbf{0},\ \mf{x}\in\partial\Omega ,\  t>0,
             \end{equation}
where $\mf{n}$ denotes the outer normal vector of $\Omega$.
The unknown function
$\rho$ is the density of the nematic liquid crystals, $\mathbf{v}$ the velocity and
$P(\rho)$ the pressure determined through the
equations of state, $\mathbf{d}$
%\in\mathbb{S}^1:=\{\mf{d}\in %\mathbb{R}^2~|~|\mf{d}|=1\}
represents the macroscopic average of
the nematic liquid crystal orientation field. The constants $\mu$, $\lambda$, $\nu$, and
$\theta$  denote the shear viscosity, the bulk viscosity, the
competition between kinetic and potential energies, and the
microscopic elastic relation time for the molecular orientation field, respectively,
they satisfy the physical conditions:
\begin{equation*} \label{0107}\mu>0,\quad \lambda +\mu\geq 0,\quad \nu>0,\quad\theta>0.\end{equation*}
 $\mathbb{I}$ denotes the $N\times N$ identity
matrix. The term $\nabla \mathbf{d}\odot\nabla \mathbf{d}$ denotes the $N\times N$ matrix
whose $(i,j)$-th entry is given by
$\partial_{x_i}\mathbf{d}\cdot\partial_{x_j}\mathbf{d}$, for $1\leq i,j\leq N$, i.e.,
 $\nabla \mathbf{d} \odot\nabla \mathbf{d}= (\nabla\mathbf{d})^{\top}\nabla\mathbf{d},$
where $(\nabla\mathbf{d})^{\top}$ denotes the transpose of the $N\times N$ matrix $\nabla\mathbf{d}$.

In 1989, Lin \cite{LFHNC} first derived a simplified
Ericksen-Leslie system modeling liquid crystal flows when the
fluid is  incompressible and viscous. Subsequently, Lin and Liu
\cite{LFHLCNC5,LFHLCPD2} established some analysis results on
the simplified Ericksen-Leslie system, such as the existence of
weak and strong solutions and the partial regularity of suitable
solutions, under the assumption that the liquid crystal  director
field is of varying length by Leslie's terminology, or variable
degree of orientation by Ericksen's terminology.

Since the supercritical nonlinearity $|\nabla\mf{d}|^2\mf{d}$
causes significant mathematical difficulties, Lin in \cite{LFHNC} introduced a
Ginzburg-Landau approximation of the simplified Ericksen-Leslie
system, i.e., $|\nabla \mf{d}|^2\mf{d}$ in (\ref{0103}) is replaced
by the Ginzburg-Landau penalty function $(1-|\mf{d}|^2)/\epsilon$ or by a
more general penalty function. Consequently, by establishing some
estimates to deal with the direction field and its coupling/interaction
with the fluid variables, a number of results on the Navier-Stokes equations
can be successfully generalized to such Ginzburg-Landau approximation model.
For examples, when $\rho$ is a constant, i.e., the homogeneous
incompressible case, Lin and Liu \cite{LFHLCNC5} proved the global
existence of weak solutions in 2D and 3D. In
particular, they also obtained the existence and uniqueness of
global classical solutions either in 2D or in 3D
for large fluid viscosity $\mu$. In addition, the existence of weak
solutions to the density-dependent incompressible flow of liquid
crystals was proved in \cite{LXZZLC3121,JFTZOGWS}. Recently, Wang
and Yu \cite{WDHYCGA}, and Liu and Qin \cite{LXGQJE} independently
established the global existence of weak solutions to
the three-dimensional compressible flow of liquid crystals
with the Ginzburg-Landau penalty function.

In the past a few years, progress has also been made on the analysis of
the model  (\ref{0101})--(\ref{0103}) by
overcoming the difficulty induced by the supercritical nonlinearity
$|\nabla\mf{d}|^2\mf{d}$. For the incompressible case,
the existence of large weak solutions in 2D was established in \cite{LFHLJYWCY} and \cite{LZLDZXY}
for a bounded domain and the whole space respectively,
and the local existence of large strong solutions and global existence of small strong solutions
in three dimensions were proved in \cite{DSJHJRHYWRZZ,HW1, LW-JDE2012,LJYDSJO,WCYWA}.
For the 3D compressible case, the existence of strong
solutions have been investigated extensively. For examples,
%Huan et al. \cite{HTWCYWHY,HTWCYWWHYSS} established
the local existence of strong solutions
and a blow-up criterion were obtained in \cite{HTWCYWHY,HTWCYWWHYSS}, while
the existence and uniqueness of global strong solutions to the Cauchy problem in
critical Besov spaces were proved in \cite{HXWHGS} provided that the initial data are close to an
equilibrium state, and the global existence of classical solutions to the Cauchy problem was shown in  \cite{LJXZHZJWG}
with smooth initial data that has small energy but possibly
large oscillations with possible vacuum and constant state as far-field condition.
Recently progress has also been made on the existence of weak solutions to multi-dimensional
problem (\ref{0101})--(\ref{0103}).
For examples, Jiang  et al \cite{JFJSWDH} established the existence of global weak solutions
to the two-dimensional problem in a bounded domain under a restriction imposed on the initial
energy including the case of small initial energy. Moreover they also obtained the existence of global
large weak solutions to the two-dimensional Cauchy problem, provided that the second component of initial
data of the direction field satisfies some geometric angle condition. At the same time, Wu and Tan \cite{WGTZGL}
established  the existence of global weak solutions to the Cauchy problem  (\ref{0101})--(\ref{0103})
by using Suen and Hoff's method \cite{SADHG}, if the initial energy around equilibrium state
is sufficiently small, the coefficients $\mu$ and
$\lambda$ satisfy $0\leq \lambda+\mu< ({3+\sqrt{21}})\mu/6$, and
the initial data $(\mf{v}_0,\mf{d}_0)$ satisfies
$\|\mathbf{v}_0\|_{\mathbf{L}^p(\mathbb{R}^3)}+\|\mathbf{d}_0\|_{\mathbf{L}^p(\mathbb{R}^3)}<\infty$
with $p>6$.

To our best knowledge, however, there are no results available on
weak solutions of the multi-dimensional problem
(\ref{0101})--(\ref{0103}) with large initial data in a bounded
domain, due to the difficulties induced by the compressibility and
the supercritical nonlinearity. It seems that the only global
existence of large weak solutions to (\ref{0101})--(\ref{0103}) was
shown in the 1D case in \cite{DSLJWCWH}. On the other hand, there
exists a  global weak solution to the multi-dimensional compressible
Navier--Stokes equations with large initial data (i.e., the initial
energy can be arbitrarily large). A question naturally arises
whether one can establish a global existence result for the problem
(\ref{0101})--(\ref{0104}) without any smallness restriction imposed
on the initial density and velocity. In the current paper, we give a
positive answer to this question in the two-dimensional case
under a restriction on the last component of initial direction field $\mf{d}_0$,
while in the three-dimensional case, a somewhat weaker existence result is obtained.

 Before stating our main result, we explain the notations and conventions
used throughout this paper. In this paper we focus our study on the case of
isentropic flows as in \cite{WDHYCGA} and assume that
\begin{equation*}
\label{0106}P(\rho)=A\rho^\gamma,\quad\mbox{ with }A>0,\ \gamma>\frac{N}{2}.
\end{equation*}
For the sake of simplicity, we define
$$I:=I_T:=(0,T),\quad Q_T=\Omega\times I,$$
\begin{equation}\label{dispassive}\mathcal{F}(t):=\mathcal{F}(\rho,\mf{v},\mf{d}):=\int_\Omega
\left(\mu |\nabla \mf{v}|^2+(\lambda+\mu)|\mm{div}\mf{v}|^2 +\theta(|\Delta
\mf{d}+|\nabla \mf{d}|^2\mf{d}|^2)\right)\mathrm{d}\mathbf{x},
\end{equation} and
\begin{equation}\label{jsw0108}
\mathcal{E}(t):= \mathcal{E}(\rho,\mf{m},\mf{d}):=\int_\Omega\left(\frac{1}{2}
\frac{|\mf{m}|^2}{\rho}1_{\{\rho>0\}}+ \frac{A}{\gamma-1}\rho^\gamma +\frac{\nu\theta|\nabla
\mf{d}|^2}{2}\right)\mm{d}\mf{x}\;\;\mbox{ with }\;\;\mf{m}=\rho\mf{v} ,
\end{equation}
where $1_{\{\rho>0\}}$ denotes the characteristic function.
We use the bold fonts to denote the product spaces, for examples,
\begin{equation*}  \begin{aligned}&
\mathbf{L}^p(\Omega):=(L^p(\Omega))^N,\quad
\mathbf{H}^k_0(\Omega):=({H}^k_0(\Omega))^N=(W^{k,2}_0(\Omega))^N,\quad
\mathbf{H}^k(\Omega):=(W^{k,2}(\Omega))^N;   \end{aligned}\end{equation*}
and  the Sobolev space with weak topology is defined as
 $$ C^0(\bar{I},\mf{L}^q_{\mm{weak}}(\Omega)):= \left\{\mf{f}: I\rightarrow
\mf{L}^q(\Omega)~\bigg|~\int_\Omega \mf{f}\cdot\mf{g}\mm{d}\mf{x}\in
C(\bar{I})\mbox{ for any }\mf{g}\in \mf{L}^\frac{q}{q-1}(\Omega)\right\}.$$
In what follows, the letter $C_0$ will denote a generic positive constant which may
depend on the dimension of space $N$, and the letter
$C(\ldots)$ will denote a generic positive constant depending on
 its variables, and is nondecreasing in its variables, except for the domain $\Omega$.
 It should be noted  that the letter $C(\ldots)$ may depend on the physical parameters
and the dimension $N$ in some places, however we usually omit this dependence for simplicity.

Our existence result of large weak solutions for (\ref{0101})--(\ref{0104}) reads as follows.
%%%%%%%%%%%%%%%%%%%%%%%%%%%%%%%%%%%%%%%%
\begin{thm}\label{thm:0101}  Let $N=2$ or $3$,
$\Omega\subset \mathbb{R}^N$ be a bounded domain of class $C^{2,\alpha}$ with $\alpha\in(0,1)$,
and the initial data $\rho_0,\,\mathbf{m}_0,\,\mathbf{d}_0$
 satisfy the following conditions:
\begin{align}   &\label{jfw0110}
  \rho_0\in L^\gamma(\Omega),\quad \rho_0\geq0\;\mbox{ a.e. in }  \Omega ,   \\
  &\label{fzu0109} \mathbf{m}_0\in \mf{L}^{\frac{2\gamma}{\gamma+1}}(\Omega),\quad
  \mf{m}_01_{\{\rho_0=0\}}=\mf{0}\;\mbox{ a.e. in }\Omega, \quad \frac{|\mathbf{m}_0|^2}{\rho_0}1_{\{\rho_0>0\}}
  \in L^1(\Omega),  \\
  &\label{jfw0112}  \mathbf{d}_0\in \mathbf{H}^2({\Omega}),\qquad
  |\mf{d}_0|\leq 1\;\mbox{ in }\;\Omega .
                                \end{align}
 Then, there exists a constant $\epsilon_0:=\epsilon_0(N,\Omega)\leq 1$ depending on $N$ and $\Omega$
 (but independent of the physical parameters in (\ref{0101})--(\ref{0103}) and the initial data),
 such that if ${d}_{0N}:=d_{0N}(\mf{x})$ (the $N$-th component of $\mf{d}_0(\mf{x})$) satisfies
\begin{equation}\label{n0109}1-d_{0N}<\epsilon_0,
\end{equation}
 the initial-boundary value problem (\ref{0101})--(\ref{0104})
has a global weak solution $(\rho,\mathbf{v},\mathbf{d})$ on $I=I_T$ for
any given $T>0$, with the following properties:
\begin{enumerate}
            \item[\quad (1)]
Regularity:
\begin{eqnarray*}%eqnarray}
&&0\leq\rho\;\mbox{ a.e. in }Q_T,\quad \rho\in C^0(\bar{I},
{L}^\gamma_{\mm{weak}}(\Omega)
  \cap C^0(\bar{I},L^p(\Omega))\cap
 L^{\gamma+\eta}(Q_T), \label{fz0114} \\
&& \mf{v}\in L^2(I,\mf{H}^1_0(\Omega)),\quad  \rho\mf{v}\in
L^\infty(I,\mf{L}^\frac{2\gamma}{\gamma+1}(\Omega))\cap
  C^0(\bar{I},\mf{L}^\frac{2\gamma}{\gamma+1}_{\mm{weak}}(\Omega)), \label{fz0115} \\
&&  |\mf{d}|\leq 1\;\mbox{ a.e. in }Q_T,\quad\mathbf{d}\in L^2(I,\mathbf{H}^2(\Omega))\cap
C^0(\bar{I},\mathbf{H}^1(\Omega)),\quad\partial_t\mathbf{d}\in {L}^{\frac{4}{3}}(I, \mathbf{L}^{2}(\Omega)),
\\
&& (\mf{n}\cdot \nabla \mathbf{d})|_{\partial\Omega}=\mathbf{0}\mbox{ in the sense of trace for a.e. $t\in I$,}
\end{eqnarray*}% eqnarray}
 where $p\in [1,\gamma)$,  and $\eta\in(0,(2\gamma-N)/N)$.
    \item[\quad (2)] Equations (\ref{0101}) and (\ref{0102}) hold in $(\mathcal{D}'(Q_T))^{N+1}$,
    and equation (\ref{0103}) holds a.e. in $Q_T$.
 \item[\quad (3)] Equation (\ref{0101}) is
satisfied in the sense of renormalized solutions, that is, $\rho$,$\mathbf{v}$ satisfy
%%%
\begin{equation}\label{0109}
\partial_tb(\rho)+\mathrm{div}[b(\rho)\mathbf{v}]+\left[\rho b'(\rho)-
b(\rho)\right]\mathrm{div}\mathbf{v}=0\mbox{ in }
\mathcal{D}'(\mathbb{R}^N\times I), \end{equation}
provided ($\rho$,$\mathbf{v}$) is prolonged to be zero on
$\mathbb{R}^N\setminus\Omega$, for any $b$ satisfying
\begin{equation*}\label{0110}  b\in C^0[0,\infty)\cap C^1(0,\infty),\quad
|b'(s)|\leq cs^{-\lambda_0},\ s\in (0,1],\ \lambda_0<1, \end{equation*}
and the growth conditions at infinity:
\begin{equation*}\label{0111}|b'(s)|\leq c\,s^{\lambda_1},\
s\geq 1,\ where\ c>0,\ 0<1+\lambda_1<\frac{(N+2)\gamma-N}{2N}.
\end{equation*}
 \item[\quad (4)] Regularity estimates:
\begin{equation}\label{jww0120fxu}\begin{aligned}
&\sup_{t\in I}(\mathcal{E}(t)+\| \mf{d}(t)-\mf{e}_N\|_{\mf{L}^2(\Omega)})
+\|(\nabla\mf{v},\nabla^2\mf{d})\|_{\mf{L}^2(Q_T)}+\|\nabla \mf{d}\|_{\mf{L}^4(Q_T)}
\leq C(\mathcal{E}_0,T,\Omega),
\\
&\|\mf{d}-\mf{e}_N\|_{\mf{L}^\infty(\Omega)}< C_0\sqrt{\epsilon_0},
\end{aligned}\end{equation}
where $\mathcal{E}_0:=\mathcal{E}(0)=\mathcal{E}(\rho_0,\mf{m}_0,\mf{d}_0)$, $\mf{e}_{2}=(0,1)$,
$\mf{e}_3=(0,0,1)$, and we have defined
$$\|(\nabla\mf{v},\nabla^2\mf{d})\|_{\mf{L}^2(Q_T)}^2=
\sum_{1\leq i,j\leq N}\|\partial_{i}v_j\|_{L^2(Q_T)}^2+\sum_{1\leq i,j,k\leq N}\|\partial_{i}
\partial_jd_k\|_{L^2(Q_T)}^2.$$
In particular, if $\Omega$ is a ball $B_R:=\{\mf{x}\in \mathbb{R}^N~|~|\mf{x}|<R\}$ with
$R\geq 1$, then
the above constant $\epsilon_0$ can be chosen to be independent of $\Omega$ for any $R\geq 1$.
Moreover, the constant $C(\mathcal{E}_0,T,\Omega)$ in \eqref{jww0120fxu}
can be replaced by a constant $C(\mathcal{E}_0,T,\|\mf{d}_0-\mf{e}_N\|_{\mf{L}^2(\Omega)})$ independent of $\Omega$.
%%%%%%%%%%%%%%%%%%%%%%%%%%%%%%
 \item[\quad (5)] In the case of $N=2$, if, in addition, $|\mf{d}_0|=1$, then the weak solution
 satisfies $|\mf{d}|\equiv 1$, and the following finite and bounded energy inequalities:
\begin{eqnarray}
&& \frac{d\mathcal{E}(t)}{dt}+\mathcal{F}(t)\leq 0 \quad\mbox{ in }\mathcal{D}'(I), \nonumber \\
&& \label{eniq1}\mathcal{E}(t)+\int_0^t\mathcal{F}(s)\mm{d}s\leq
\mathcal{E}_0 \quad\mbox{ for a.e. }t\in I.
\end{eqnarray}
 \end{enumerate}\end{thm}

 \iffalse
\begin{rem}
In Theorem \ref{thm:0101}, we have used notations on Sobolev spaces that
\begin{equation*}  \begin{aligned}&
\mathbf{H}^k_0(\Omega):=({H}^k_0(\Omega))^2=(W^{k,2}_0(\Omega))^2,\quad
\mathbf{L}^p(\Omega):=(L^p(\Omega))^2\mbox{ for }p=2,\ {2\gamma}/({\gamma+1}),\\
&\mf{L}^\frac{2\gamma}{\gamma+1}_{\mm{weak}}(\Omega):
=({L}^\frac{2\gamma}{\gamma+1}_{\mm{weak}}(\Omega))^2,\quad
\mathbf{H}^k(\Omega):=(W^{k,2}(\Omega))^2.\end{aligned}\end{equation*}
%%%
Other notations with bold form will appear in the rest of this paper
(to denote the product spaces). For the convenience of readers, we
give the definition of the Sobolev spaces with weak topology: $
C^0(\bar{I},\mf{L}^q_{\mm{weak}}(\Omega)):= \{\mf{f}: I\rightarrow
\mf{L}^q(\Omega)~|~\int_\Omega \mf{f}\cdot\mf{g}\mm{d}\mf{x}\in
C(\bar{I})\mbox{ for any }\mf{g}\in \mf{L}^\frac{q}{q-1}(\Omega)\}$.
\end{rem}
\fi

\begin{rem}
The proof of Theorem \ref{0101} remains basically unchanged
if the motion of the fluid is driven by a bounded external force,
i.e., when the momentum equations (\ref{0102}) contain an additional term
$\rho\mf{f}(\mf{x},t)$ with $\mf{f}$ being a bounded and measurable function.
We remark here that we do not require any smallness condition on $\mf{f}$.
However, we are not clear whether the above theorem still holds with
non-homogenous boundary condition in place of Neumann boundary condition
``$(\mf{n}\cdot\nabla \mf{d})|_{\partial\Omega}=\mathbf{0}$''.
In the proof of Theorem \ref{thm:0101}, we use the Neumann boundary condition only in order to
deduce a maximum principle on $\mf{d}$.
\end{rem}
%%%%%%%%%%%%%%%%%%%%%%%%%
\begin{rem}
We mention that the regularity requirement ``$\mathbf{d}_0(\mathbf{x})\in \mathbf{H}^2({\Omega})$''
is not optimal, for example, if we have
``$\mathbf{d}_0(\mathbf{x})\in \mathbf{W}^{1,p}({\Omega})$ with $p>N$", then the above theorem still
holds, and this can be shown by a standard approximate approach. On the other hand,
we do not known whether ``$\mathbf{d}_0(\mathbf{x})\in \mathbf{H}^{1}({\Omega})$'' is
the lowest regularity requirement, since it involves the problem of the Sobolev maps between two manifolds
with the lower boundedness condition (\ref{n0109}).
\end{rem}
\begin{rem} In view of the above regularity estimates in a ball, we
can make use of a domain expansion technique to obtain a similar existence
result of global weak solutions to the corresponding Cauchy problem, for
which the expression of energy $\mathcal{E}(t)$ should be written in
a form around some equilibrium state
$(\rho_\infty,\mf{v}_\infty,\mf{e}_N)$ with $\rho_\infty>0$ to make
the energy integral sense (see \cite[Theorem 1.2]{JFJSWDH}). Of course in this case,
we can also establish a similar existence result of global weak
solutions to the corresponding incompressible problem.
\end{rem}

We now describe the main idea of the proof of Theorem
\ref{thm:0101}. For the Ginzburg-Landau approximation model to
(\ref{0101})--(\ref{0103}), based on some new estimates to deal with
the direction field and its coupling/interaction with the fluid variables,
Wang and Yu in \cite{WDHYCGA} adopted a classical
three-level approximation scheme which consists of the
Faedo-Galerkin approximation, artificial viscosity, an artificial
pressure and the celebrated weak continuity of the effective
viscous flux to overcome the difficulty of possible large
oscillations of the density, and established the existence of weak
solutions. These techniques were developed in \cite{LPLMTFM98} and \cite{FENAPHOFJ35801,JSZPO}
for the compressible Navier-Stokes equations, we refer to the monograph
\cite{NASII04} for more details. In the proof of Theorem \ref{thm:0101},
we also adopt the three-level approximation scheme, so the key steps are to deduce the
\emph{a priori} estimates and to construct approximate solutions to the third approximate problem.
Compared with the Ginzburg-Landau approximation model in
\cite{WDHYCGA}, however, the system (\ref{0101})--(\ref{0103}) is
much more difficult to deal with, due to the supercritical nonlinearity $|\nabla
\mf{d}|^2\mf{d}$ in (\ref{0103}). Consequently, not like that in \cite{WDHYCGA},
one can not deduce the (sufficiently) strong estimate $\nabla^2\mf{d}\in L^2(I,\mf{L}^2(\Omega))$
directly from the basic energy inequality (\ref{eniq1}).
Recently, Ding and Wen obtained the global existence
and uniqueness of strong solutions to the 2D density-dependent
incompressible model with small initial energy and positive initial
density away from zero in \cite{DSJWHYSN}, where they got
$\nabla^2\mf{d}\in L^2(I,\mf{L}^2(\Omega))$ from the basic energy
inequality under the smallness condition of the initial energy.
In fact, they first deduced ${\nu\|\nabla\mf{d}\|_{\mf{L}^2(\Omega)}^2}/{2}\leq \mathcal{E}_0$ and
$\nu\theta\|\Delta \mf{d}\|_{\mf{L}^2(Q_T)}^2\leq
\mathcal{E}_0+\nu\theta\||\nabla \mf{d}||_{{L}^4(Q_T)}^4$ by employing
the basic energy inequality, and then made use of the inequality
\begin{equation*}\label{jww0121}
\|\nabla \mf{d}\|_{\mf{L}^4(\Omega)}^4\leq {C}(\Omega)(\|\nabla^2 \mf{d}\|_{\mf{L}^2(\Omega)}^2\|\nabla\mf{d}\|_{\mf{L}^2(\Omega)}^2+\|\nabla \mf{d}\|_{\mf{L}^{2}
(\Omega)}^4)
\end{equation*}  for some constant $C(\Omega)$ depending on $\Omega\subset \mathbb{R}^2$,
which follows from the elliptic estimates and an interpolation inequality (see \cite[Lemma 2.4]{DSJWHYSN}),
to infer that
\begin{equation}\label{jww0122}\|\nabla^2\mf{d}
\|_{\mf{L}^2(Q_T)}^2\leq C(\mathcal{E}_0,\|\nabla^2\mf{d}_0\|_{\mathbf{L}^2(\Omega)}), \end{equation}
provided that the initial energy is sufficiently small. Motivated by this study,
Jiang et al \cite{JFJSWDH} established the global existence of weak solutions to the corresponding compressible problem.
In the current paper, we shall use another version of the interpolation inequality
\begin{equation*} %\label{}
\|\nabla \mf{d}\|_{\mf{L}^4(\Omega)}^4\leq {C}(\Omega)(\|\nabla^2 \mf{d}\|_{\mf{L}^2(\Omega)}^{{2}}
\|\mf{d}-\mf{e}_N\|_{\mf{L}^\infty(\Omega)}^2+\| \mf{d}-\mf{e}_N\|_{\mf{L}^{4} (\Omega)}^{4})
\end{equation*} for some constant $C(\Omega)$ depending on $\Omega\subset \mathbb{R}^N$ ($N=2$ or $3$),
from which the estimate (\ref{jww0122}) can also be deduced  for $\Omega\subset \mathbb{R}^N$ if $\|\mf{d}-\mf{e}_N\|_{\mf{L}^\infty(\Omega)}^2$
is sufficiently small. Now the question is whether the smallness of
$\|\mf{d}-\mf{e}_N\|_{\mf{L}^\infty(\Omega)}^2$ is guaranteed by smallness of the initial
data $\|\mf{d}_0-\mf{e}_N\|_{\mf{L}^\infty(\Omega)}^2$.
Fortunately, this is the case by applying the maximum principle to nonnegative lower bounds of solutions to the equations (\ref{0103})
and the condition $|\mf{d}|\leq 1$.
Consequently, we deduce the desired energy estimates on $\mf{d}$
 from the energy inequality. %% as in \cite{WDHYCGA}.
With these estimates in hand, we can adopt and modify the three-dimensional approximation scheme approach
to show Theorem \ref{thm:0101}, if we can construct a solution to the following third approximate problem:
\begin{eqnarray}&&%\label{}
\partial_t\rho+\mathrm{div}(\rho\mathbf{v})
=\varepsilon\Delta\rho ,   \\
&& \label{111111sd}\partial_t\mathbf{d}+\mathbf{v}\cdot \nabla \mathbf{d}
= \theta(\Delta \mathbf{d}+f_{\varepsilon}(|\nabla \mathbf{d}|^2)\mf{d}),  \\[1mm]
&& \begin{aligned}%\label{}
&\int_\Omega
(\rho\mathbf{v})(t)\cdot\mathbf{\mathbf{\Psi}}\mathrm{d}\mathbf{x}
-\int_{\Omega}\mathbf{m}_0\cdot\mathbf{\mathbf{\Psi}}\mathrm{d}\mathbf{x} \\
& = \int_0^t\int_\Omega\bigg[\mu\Delta\mathbf{v}+(\mu+\lambda)\nabla
\mm{div}\mf{v}-\nabla P -\delta \nabla
\rho^\beta-\varepsilon(\nabla\rho\cdot\nabla \mathbf{v})  \\
&\qquad\qquad -\mathrm{div}(\rho\mathbf{v}\otimes\mathbf{v})
-\nu\mathrm{div}\left(\nabla \mathbf{d}\otimes\nabla\mathbf{d} -\frac{|\nabla\mf{d}|^2
\mathbb{I}}{2}\right)\bigg]\cdot\mathbf{\Psi}\mathrm{d}\mathbf{x}\mathrm{d}s,
 \end{aligned} \end{eqnarray}
 where the $n$-dimensional Euclidean space $X_n$ will be introduced in Section \ref{sec:jww03}, $\varepsilon$, $\delta$, $\beta> 0$ are
constants, and the smooth function $f_\varepsilon(x)\geq 0$  satisfying
\begin{equation}\label{appfun}f_{\varepsilon}(x)=x\mbox{ if }N=2;\qquad  \left\{
                  \begin{array}{l}
                    0\leq f'(x)\leq 1,  \\[0.5em]
                     f(x)=\varepsilon^{-1}\mbox{ if }x\geq \varepsilon^{-1}, \\[0.5em]
                    0\leq x-f_\varepsilon(x)\rightarrow 0\mbox{ as }\varepsilon\rightarrow 0,
                  \end{array}
                \right.\mbox{ if }N=3.
\end{equation}

It should be noted that  the third approximate
problem above still enjoys the desired energy estimates (see Proposition \ref{pro:energy}), thus it is easy to establish
the unique solvability of the third approximate problem  in the 2D case by following the same proof
as in \cite{JFJSWDH}. However, the proof in \cite{JFJSWDH} can not be directly applied to the 3D case, and
the difficulty lies in that we could not deduce a global estimate on
$\|\partial_t\mf{d}\|_{L^\infty(I,\mf{L}^2(\Omega))}$ (see (\ref{nnn0457}))
for the 3D approximate problem (\ref{n301})--(\ref{n0306}). To overcome this difficulty,
we introduce the cut-off function \eqref{appfun} to get a global estimate of
$\|\partial_t\mf{d}\|_{L^\infty(I,\mf{L}^2(\Omega))}$. On the other hand, we have to pay the price for this, namely,
for the 3D approximate problem \eqref{111111sd} based on a cut-off function with Neumann boundary condition, we can not
show $|\mf{d}|=1$ when $|\mf{d}_0|=1$. This is the reason why the
solution in Theorem \ref{thm:0101} does not satisfy $|\mf{d}|=1$ in three dimensions.
%% Maybe we need to explore a new approximate method and technique of proof.

 The rest of paper is organized as follows. In Section \ref{sec:02} we deduce the basic energy equalities
 from the third approximate problem and derive more energy estimates on
$\mf{d}$ under the assumption (\ref{n0109}). In Section \ref{subsystem} we introduce the strong solvability
of sub-systems in the third approximate problem, while the unique
solvability of the third approximate problem is established in Section \ref{sec:jww03}.
Finally, we briefly sketch how to use the standard three-level approximation
scheme to prove Theorem \ref{thm:0101} in Section \ref{sec:03}.

\section{\emph{A priori} for the third approximate problem}\label{sec:02}
This section is devoted to formal derivation of the \emph{a priori} energy estimates for
the third approximate equations:
\begin{eqnarray}
&& \label{jww0101}\partial_t\rho+\mathrm{div}(\rho\mathbf{v})=\varepsilon\Delta\rho,\\
 &&  \label{jww0102} \partial_t
(\rho\mathbf{v})+\mm{div}(\rho\mathbf{v}\otimes \mathbf{v})+ \nabla
P+\delta \nabla
\rho^\beta+\varepsilon(\nabla\rho\cdot\nabla \mathbf{v})\nonumber\\
&&\quad  =\mu \Delta\mathbf{v}+(\mu+\lambda)\nabla
\mm{div}\mf{v}
-\nu\mm{div}\Big( \nabla \mathbf{d}\odot\nabla
\mathbf{d}-\frac{1}{2}|\nabla \mf{d}|^2\mathbb{I}\Big),  \\
&&  \label{jww0103}\partial_t\mathbf{d}+\mathbf{v}\cdot \nabla
\mathbf{d}= \theta(\Delta \mathbf{d}+f_\varepsilon(|\nabla \mathbf{d}|^2)\mf{d}),
\end{eqnarray}
in a bounded domain $\Omega\subset \mathbb{R}^N$ with initial data
\begin{eqnarray}\label{jwwn0305}&&\rho(\mathbf{x},0)=\rho_0>0,
\quad\mathbf{d}(\mathbf{x},0)=\mathbf{d}_0,\quad
 \mathbf{v}(\mathbf{x},0)=\mf{v}_0,\end{eqnarray}
 and boundary conditions
 \begin{eqnarray}\label{jwwn0306}
 \nabla\rho\cdot\mf{n}|_{\partial\Omega}=0,\quad
\mathbf{v}|_{\partial \Omega}=\mathbf{0},\quad
(\mf{n}\cdot\nabla\mathbf{d})|_{\partial\Omega}=\mathbf{0}.\end{eqnarray}
The\emph{ a priori }estimates will play a crucial role in the proof of existence.
We consider a classical solution $(\rho, \mathbf{v}, \mathbf{d})$ of
the initial-boundary problem (\ref{jww0101})--(\ref{jwwn0306}) with $\rho>0$.

\subsection{Basic energy estimates}
We first deduce some basic energy estimates without any smallness condition imposed on the initial data.
%%%%%%%%%%%%%%%%%%%%%%%%%%%%%%%%%%%%%%%%%%%%%%%%%%%%%%%%%%%
\subsubsection{Maximum principle on $|\mf{d}|$}
The macroscopic average of the nematic liquid crystal orientation field $\mf{d}$ satisfies
\begin{equation} \label{0201}
|\mf{d}|\leq 1\mbox{ in }Q_T, \mbox{ if }|\mf{d}_0|\leq 1\mbox{ in }\Omega\subset \mathbb{R}^N.
\end{equation}
Next, we give a proof of (\ref{0201}) for the reader's convenience.
Multiplying the $\mf{d}$-system (\ref{0103}) by $\mf{d}$, we obtain
\begin{equation*} \label{0202}\begin{aligned}
\frac{1}{2}\partial_t|\mf{d}|^2+\frac{1}{2}\mf{v}\cdot\nabla
|\mf{d}|^2=\theta(\Delta\mf{d}\cdot\mf{d}+f_\varepsilon(|\nabla
\mf{d}|^2)|\mf{d}|^2).\end{aligned}\end{equation*}
From the identity $\Delta|\mf{d}|^2=2|\nabla\mf{d}|^2+2\Delta\mf{d}\cdot\mf{d}$ it follows that
\begin{equation} \label{0205}\begin{aligned}
\partial_t(|\mf{d}|^2-1)-\theta\Delta(|\mf{d}|^2-1)\leq -\mf{v}\cdot\nabla
(|\mf{d}|^2-1)+2\theta |\nabla\mf{d}|^2(|\mf{d}|^2-1).\end{aligned}
\end{equation}
Now, letting $d=|\mf{d}|^2-1$ and $d_+=\max\{d,0\}\geq 0$, multiplying  (\ref{0205}) by $d_+$
and integrating over $\Omega$, we integrate by parts and use the boundary conditions to infer that
$$ \frac{d}{dt}\int_\Omega d_+^2\mm{d}\mf{x}\leq \int_\Omega(4\theta|\nabla \mf{d}|^2
+\mm{div}\mf{v})d_+^2\mm{d}\mf{x} \leq \|4\theta|\nabla \mf{d}|^2+|\mm{div}\mf{v}|\|_{L^\infty(\Omega)}
\int_\Omega d_+^2\mm{d}\mf{x}. $$
Assuming that $(\mf{v},\mf{d})$ satisfies the following regularity
\begin{equation*}
\|4\theta|\nabla
\mf{d}|^2+|\mm{div}\mf{v}|\|_{L^1(I,L^\infty(\Omega))}<\infty,
\end{equation*}
we are able to apply Gronwall's inequality to get (\ref{0201}) immediately.

We shall see that all the couples
$(\mf{v}_n,\mf{d}_n)$ in the third approximate solutions constructed in Section \ref{sec:jww03} satisfy the regularity
required above. We remark that (\ref{0205}) becomes an equality in the two-dimensional case, and one
can get by directly multiplying (\ref{0205}) with $|\mf{d}|^2-1$ that
\begin{equation}|\mf{d}|=1\mbox{ in }Q_T, \quad\mbox{ if }|\mf{d}_0|=1\mbox{ in }\Omega\subset \mathbb{R}^2.
\end{equation}

\subsubsection{Energy inequality}
 Integrating by parts and utilizing the boundary conditions, one easily sees that the system
 (\ref{jww0101})--(\ref{jww0102}) satisfies the energy conservation:
\begin{equation}\begin{aligned}\label{basicenergy}
& \frac{d}{dt}{{E}_\delta}(t) +\int_\Omega \left(\mu |\nabla
\mf{v}|^2+(\lambda+\mu)|\mm{div}\mf{v}|^2
+\varepsilon\delta\beta\rho^{\beta-2}|\nabla \rho|^2
+{A\varepsilon\gamma}\rho^{\gamma-2}|\nabla\rho|^2\right)\mathrm{d}\mathbf{x}\\
&=-\nu\int_\Omega (\nabla
\mf{d})^T\Delta \mf{d}\cdot\mf{v}\mathrm{d}\mathbf{x},
\end{aligned}\end{equation}
where $(\nabla\mf{d})^T\Delta \mf{d}:=(\partial_id_j)_{N\times N}\Delta \mf{d}$ and
$$\displaystyle{{E}}_\delta(t)
=\int_{\Omega}\left(\frac{1}{2}\frac{|\mathbf{m}|^2}{\rho}1_{\{\rho>0\}}+\frac{A\rho^\gamma}{\gamma-1}
+\frac{\delta}{\beta-1}\rho^\beta\right)\mm{d}\mf{x}\;\;\mbox{ with }\;\;\mf{m}=\rho\mf{v}.$$

 Multiplying (\ref{jww0103}) by $-\Delta \mf{d}$ % in $L^2(\Omega )$
  and integrating by parts, we have
 \begin{equation*}\begin{aligned}%\label{}
 \frac{1}{2}&\frac{d}{dt}
\int_\Omega|\nabla \mf{d}|^2\mm{d}\mf{x}+\theta\int_\Omega |\Delta
\mf{d}|^2\mm{d}\mf{x}\\
&=\int_\Omega(\mf{v}\cdot\nabla
\mf{d})\cdot\Delta\mf{d}\mm{d}\mf{x}+\theta\int_\Omega \left(f_\varepsilon(|\nabla
\mf{d}|^2)|\nabla\mf{d}|^2+2\sum_{1\leq i,j\leq N}f_\varepsilon'(|\nabla
\mf{d}|^2)(\partial_{x_i} {d}_j\nabla\partial_{x_i}{d}_j\cdot\nabla)\mf{d}\cdot\mf{d}\right)\mm{d}\mf{x}\\
&\leq \int_\Omega(\mf{v}\cdot\nabla
\mf{d})\cdot\Delta\mf{d}\mm{d}\mf{x}+\theta\int_\Omega  (|\nabla
\mf{d}|^4+C_0|\nabla\mf{d}|^2|\nabla^2\mf{d}|)\mm{d}\mf{x},\end{aligned}\end{equation*}
which, together with \eqref{basicenergy}, implies
\begin{equation}\label{0542}\begin{aligned}&\frac{d}{dt}{\mathcal{E}}_\delta(t)+\|\mu |\nabla
\mf{v}|^2+(\lambda+\mu)|\mm{div}\mf{v}|^2  +\nu\theta|\Delta \mf{d}|^2+\varepsilon\delta\beta\rho^{\beta-2}|\nabla
\rho|^2+{A\varepsilon\gamma}\rho^{\gamma-2}|\nabla\rho|^2\|_{L^1(\Omega)}\\
& \leq \theta\nu\int_\Omega  (|\nabla
\mf{d}|^4+C_0|\nabla\mf{d}|^2|\nabla^2\mf{d}|)\mm{d}\mf{x},
\end{aligned}\end{equation}
where
$$\displaystyle {\mathcal{E}}_\delta(t):={\mathcal{E}}_\delta(\rho,\mf{m},\mf{d}):
=\int_{\Omega}\left(\frac{1}{2}\frac{|\mathbf{m}|^2}{\rho}1_{\{\rho>0\}}+\frac{A}{\gamma-1}\rho^\gamma
+\frac{\delta}{\beta-1}\rho^\beta+\frac{\nu|\nabla
\mf{d}|^2}{2}\right)\mm{d}\mf{x}.$$

For the two-dimensional case, recalling $|\mf{d}|\equiv 1$ and $f_\varepsilon({x})\equiv {x}$ for $x\geq 0$,
we can deduce the standard energy equality.
In fact, multiplying the equations \eqref{jww0103} by $\Delta\mf{d}+|\nabla\mf{d}|^2\mf{d}$,
integrating by parts and using the boundary conditions, one deduces that
\begin{equation*}%\label{}
\frac{1}{2}\frac{d}{dt}
\int_\Omega|\nabla \mf{d}|^2\mm{d}\mf{x}+\theta\int_\Omega |\Delta \mf{d}+|\nabla
\mf{d}|^2\mf{d}|^2\mm{d}\mf{x}=\int_\Omega(\mf{v}\cdot\nabla
\mf{d})\cdot\Delta\mf{d}\mm{d}\mf{x},\end{equation*}
which, together with \eqref{basicenergy}, yields
\begin{equation}\label{05420}\begin{aligned}&\int_\Omega [\mu |\nabla \mf{v}|^2
+(\lambda+\mu)|\mm{div}\mf{v}|^2  +\nu\theta(|\Delta \mf{d} +|\nabla\mf{d}|^2\mf{d}|^2)\\
&\qquad +\varepsilon\delta\beta\rho^{\beta-2}|\nabla \rho|^2
+{A\varepsilon\gamma}\rho^{\gamma-2}|\nabla\rho|^2]\mathrm{d}\mathbf{x}+\frac{d}{dt}{\mathcal{E}}_\delta(t)= 0,
\end{aligned}\end{equation}
whence,
 \begin{equation}\label{05421}
 \begin{aligned}&\int_0^t\int_\Omega [\mu |\nabla \mf{v}|^2
+(\lambda+\mu)|\mm{div}\mf{v}|^2  +\nu\theta(|\Delta \mf{d} +|\nabla\mf{d}|^2\mf{d}|^2) \\
&\qquad\quad  +\varepsilon\delta\beta\rho^{\beta-2}|\nabla \rho|^2
+{A\varepsilon\gamma}\rho^{\gamma-2}|\nabla\rho|^2]\mathrm{d}\mathbf{x}\mm{d}s
+{\mathcal{E}}_\delta(t) = {\mathcal{E}}_{\delta,0},
\end{aligned}\end{equation}
where ${\mathcal{E}}_{\delta,0}:={\mathcal{E}}_\delta(\rho_0,\mf{m}_0,\mf{d}_0)$.

\subsubsection{Maximum principle on lower bounds}
 The system (\ref{jww0103}) possesses the following maximum principle on nonnegative lower bounds:
\begin{equation}\label{jww0205}\begin{aligned}
d_i\geq \underline{d}_{0i}\;\;\mbox{ if }\;\; {d}_{0i}\geq \underline{d}_{0i}\geq 0
\;\;\mbox{ for any given constant }\;\;\underline{d}_{0i},
\end{aligned}\end{equation}
where $1\leq i \leq N$, and we have denoted the $i$-th component of $\mf{d}$ and $\mf{d}_0$ by ${d}_i$
and $d_{0i}$, respectively. This conclusion will play a crucial role in this paper, so we give
its proof here for the reader's convenience.

Letting
$$\omega_i=d_i-\underline{d}_{0i}, \quad \omega^-_i=\min\{\omega_i,0\}\leq 0,$$
 we can deduce from (\ref{jww0103}) that
\begin{equation}\label{jww0206}
\partial_t\omega_i-\theta \Delta\omega_i=\theta f_\varepsilon( |\nabla
\mf{d}|^2)(\omega_i+\underline{d}_{0i})-\mf{v}\cdot\nabla \omega_i.
\end{equation}
Multiplying (\ref{jww0206}) by $\omega^-_i$%in $L^2(\Omega)$
, and using the Neumann boundary condition,  (\ref{appfun})
and H\"older's inequality, we find that
\begin{equation*}\begin{aligned}
& \frac{1}{2}\frac{d}{dt}\|\omega^-_i\|_{L^2(\Omega)}^2 +\theta \|\nabla
\omega^-_i\|_{L^2(\Omega)}^2  \\
& = \int_{\Omega} [\theta f_\varepsilon(|\nabla\mf{d}|^2)(\omega_i+\underline{d}_{0i})
-\mf{v}\cdot\nabla\omega_i]\omega^-_i\mm{d} \mf{x}\\
&= \int_{\Omega}\theta f_\varepsilon( |\nabla\mf{d}|^2|)|\omega_i^-|^2\mm{d} \mf{x}
+\frac{1}{2}\int_{\Omega}\mm{div}\mf{v}|\omega_i^-|^2\mm{d} \mf{x}
+\int_{\Omega} \theta f_\varepsilon( |\nabla\mf{d}|^2)\underline{d}_{0i}\omega^-_i\mm{d} \mf{x},  \\
& \leq \Big( \left\|\theta |\nabla
\mf{d}|^2\right\|_{L^\infty(\Omega)}+\frac{1}{2}\|\mm{div}\mf{v}\|_{L^\infty(\Omega)}\Big)
\|\omega^-_i\|_{L^2(\Omega)}^2+\int_{\Omega} \theta |\nabla\mf{d}|^2\underline{d}_{0i}
\omega^-_i\mm{d} \mf{x},
\end{aligned} \end{equation*}
which, together with the fact $\underline{d}_{0i}
\omega^-_i\leq 0$, yields
\begin{equation*}\begin{aligned}\frac{d}{dt}\|\omega^-_i
\|_{L^2(\Omega)}^2 \leq \left(2\left\|\theta |\nabla
\mf{d}|^2\right\|_{L^\infty(\Omega)}+\|\mm{div}\mf{v}\|_{L^\infty(\Omega)}\right)
\|\omega^-_i\|_{L^2(\Omega)}^2 .
\end{aligned}\end{equation*}
Hence, if we apply Gronwall's inequality to the above inequality, we obtain
\begin{equation*}\begin{aligned}
\|\omega^-_i(t)\|_{L^2(\Omega)}^2 \leq \|\omega^-_i(0) \|_{L^2(\Omega)}^2 e^{{\int_0^t(2\left\|\theta |\nabla
\mf{d}|^2\right\|_{L^\infty(\Omega)}+\|\mm{div}\mf{v}\|_{L^\infty(\Omega)})\mm{d}s}}=0,
\end{aligned}\end{equation*}
which gives (\ref{jww0205}).

\subsection{More estimates under the small oscillation condition imposed on $\mf{d}$}\label{sec:jww0202}
To obtain more estimates on $\mf{d}$ under the small oscillation condition, we first introduce
the well-known Nirenberg interpolation inequality (see \cite[Theorem]{Nirenberg}):
\begin{lem}\label{Nirenberg}
Let $u$ belong to $L^q(\mathbb{R}^N)$ and its derivatives of order $m$, $\nabla^m u$, belong to
 $L^r(\mathbb{R}^N)$, $1\leq q$, $r\leq \infty$. Then for the derivatives $\nabla^ju$, $0\leq j<m$,
 the following inequality holds.
\begin{equation}\label{Nirenbergin}\|\nabla^ju\|_{L^p(\mathbb{R}^N)}\leq C_0\|\nabla^mu\|_{L^r(\mathbb{R}^N)}^\alpha
 \|u\|_{L^q(\mathbb{R}^N)}^{1-\alpha},\end{equation}
 where
 $$\frac{1}{p}=\frac{j}{n}+\alpha\left(\frac{1}{r}-\frac{m}{n} \right)+(1-\alpha)\frac{1}{q},$$
for all $\alpha$ in the interval
  $$\frac{j}{m}\leq \alpha \leq 1$$
(the constant $C_0$ depends only on $n$, $m$, $j$, $q$, $r$, $\alpha$), with the following exceptional cases:
\begin{enumerate}[\quad (1)]
  \item If $j=0$, $rm<n$ and $q=\infty$, then we make the additional assumption that either $u$ tends to zero
  at infinity or $u\in L^{\tilde{q}}(\mathbb{R}^N)$ for some finite $\tilde{q}>0$.
  \item If $1<r<\infty$, and $m-j-n/r$ is a non-negative integer,
then (\ref{Nirenbergin}) holds only for $\alpha$ satisfying $j/m\leq \alpha<1$.
\end{enumerate}

In addition, for a bounded domain $\Omega$ (with smooth boundary) the above assertions hold
if we add to the right side (\ref{Nirenbergin}) the term
$$C(\Omega)\|u\|_{L^{\tilde{q}}(\Omega)}$$
for any $\tilde{q}\geq 1$. All the relevant constants thus depend also on the domain.
\end{lem}

Next, we derive more estimates on $\mf{d}$ % in the three-dimensional case
under
the assumption that the initial value of $d_N$ satisfies
\begin{equation}\label{jww0212}
1-\epsilon_0\leq  d_{0N}\leq 1,\quad\;\mbox{ for some }\;\; \epsilon_0\in (0,1]. \end{equation}
It should be noted that the constant $C(\Omega)$ in the following deduction will denote various positive constants depending on
 its variable $\Omega$, but the constants $\tilde{C}_0$ and $\tilde{C}_1(\Omega)$--$\tilde{C}_3(\Omega)$ are fixed.

First, one gets from the maximum principle that
\begin{equation}\label{jww0213}\begin{aligned}
1-\epsilon_0\leq  d_{N}(\mbox{x},t)\leq 1\;\;\mbox{ for any }\; t>0\; \mbox{ and any }\; \mf{x}\in\Omega .
\end{aligned}\end{equation}
 Recalling $|{\mf{d}}|\leq  1$, that is $\sum_{i=1}^N{d}_i^2 \leq 1$, one finds that
%%% \begin{equation}\label{jww0214}\begin{aligned}
$$ |{d}_i|\leq \sqrt{1-d_N^2}<\sqrt{2\epsilon_0}
\quad\mbox{ for }1\leq i\leq N-1, $$
which combined with (\ref{jww0213}) leads to
\begin{equation}\label{jww0216bounded}\begin{aligned}
\|\mf{d}-{\mf{e}}_N\|_{\mathbf{L}^\infty(\Omega)}\leq \tilde{C}_0\sqrt{\epsilon_0}.
\end{aligned}\end{equation}

Thanks to Lemma \ref{Nirenberg}, we have
\begin{eqnarray}\label{jww0208}
&&\|\nabla {u}\|_{L^4(\Omega)}\leq
C(\Omega)(\| \nabla^2 u\|_{{L}^{2}(\Omega)}^{\frac{1}{2}}\|u\|_{L^\infty(\Omega)}^{\frac{1}{2}}
+\|u\|_{L^4(\Omega)})\mbox{ for }\Omega\subset \mathbb{R}^N\mbox{ with }N=2\mbox{ or }3,\nonumber \end{eqnarray}
which yields
\begin{equation}  \label{jww0216}
\|\nabla \mf{d}\|_{\mathbf{L}^4(\Omega)}^4\leq
C(\Omega)(\| \nabla^2 \mf{d}\|_{\mathbf{{L}}^{2}(\Omega)}^{2}\|\mf{d}-{\mf{e}}_N\|_{\mathbf{L}^\infty(\Omega)}^{2}
+\|\mf{d}-{\mf{e}}_N\|_{\mathbf{L}^4(\Omega)}^{4}).  \end{equation}

To bound the right hand of (\ref{jww0216}), we  shall use
the following elliptic estimate: There exists a constant $\tilde{C}_1(\Omega)$, such that
 \begin{equation}\begin{aligned}\label{jww0218}
 \|\nabla^2 \mf{d}\|_{\mathbf{{L}}^2(\Omega)}\leq \tilde{C}_1(\Omega)(\|\Delta
 \mf{d}\|_{\mathbf{{L}}^2(\Omega)}+\|\nabla\mf{d}\|_{\mathbf{{L}}^2(\Omega)})
 \quad\mbox{for any }\nabla \mf{d}\in \mf{H}^1(\Omega)  \end{aligned}\end{equation}
where $\Omega\subset\mathbb{R}^N$, which can be deduced from \cite[Lemma 4.27]{NASII04}.
 Thus, putting (\ref{jww0216bounded})--(\ref{jww0218}) together, we conclude that
\begin{equation}\begin{aligned}\label{jww0219}
\|\nabla\mf{d}\|_{\mathbf{L}^4(\Omega)}^4\leq &
\tilde{C}_2(\Omega)\tilde{C}_0\epsilon_0\| \Delta \mf{d}\|_{\mathbf{{L}}^{2}(\Omega)}^{2}+C(\Omega)(
\|\mf{d}-{\mf{e}}_N\|_{\mathbf{L}^\infty(\Omega)}^{2}\|\nabla\mf{d}\|_{\mathbf{{L}}^2(\Omega)}^2+
\|\mf{d}-{\mf{e}}_N\|_{\mathbf{L}^4(\Omega)}^{4})\\
\leq & \tilde{C}_2(\Omega)\tilde{C}_0\epsilon_0\| \Delta \mf{d}\|_{\mathbf{{L}}^{2}(\Omega)}^{2}
+C(\Omega)(\|\nabla\mf{d}\|_{\mathbf{{L}}^2(\Omega)}^2+ \|\mf{d}-{\mf{e}}_N\|_{\mathbf{L}^2(\Omega)}^{2}),
 \end{aligned}\end{equation}
where the constant $\tilde{C}_2(\Omega)\geq 2^{-1}$ only depends on $\Omega$.
Utilizing (\ref{jww0219}), (\ref{jww0218}), and Cauchy's and H\"older's inequalities,
we can deduce from (\ref{0542}) that
 \begin{equation}\begin{aligned}\label{energyequality}
&\left\|\mu |\nabla \mf{v}|^2+(\lambda+\mu)|\mm{div}\mf{v}|^2  +\nu\theta|\Delta \mf{d}|^2
+\varepsilon\delta\beta\rho^{\beta-2}|\nabla
\rho|^2+{A\varepsilon\gamma}\rho^{\gamma-2}|\nabla\rho|^2\right\|_{L^1(\Omega)}
+\frac{d}{dt}{\mathcal{E}}_\delta(t) \\
& \leq  \tilde{C}_3\theta\nu(\|\nabla
\mf{d}\|^4_{\mf{L}^4(\Omega )} +\|\nabla
\mf{d}\|^2_{\mf{L}^4(\Omega )}\|\nabla^2\mf{d}\|_{\mf{L}^2(\Omega)})\\
&\leq \tilde{C}_3\theta\nu\left[\tilde{C}_2(\Omega)\tilde{C}_0\epsilon_0
\|\Delta\mf{d}\|_{\mathbf{{L}}^{2}(\Omega)}^{2}+C(\Omega)(\|\nabla\mf{d}\|_{\mathbf{{L}}^2(\Omega)}^2+
\|\mf{d}-{\mf{e}}_N\|_{\mathbf{L}^2(\Omega)}^{2})\right.\\
&\qquad\qquad  +\tilde{C}_1(\Omega)(\sqrt{\tilde{C}_2(\Omega)\tilde{C}_0\epsilon_0}
\|\Delta \mf{d}\|_{\mathbf{{L}}^{2}(\Omega)}+\sqrt{C(\Omega)}(\|\nabla\mf{d}\|_{\mathbf{{L}}^2(\Omega)}+
\|\mf{d}-{\mf{e}}_N\|_{\mathbf{L}^2(\Omega)})\\
&\qquad \qquad\quad \left. \times (\|\Delta
 \mf{d}\|_{\mathbf{{L}}^2(\Omega)}+\|\nabla\mf{d}\|_{\mathbf{{L}}^2(\Omega)})\right]  \\
&\leq  \tilde{C}_3\theta\nu\left[\left(\tilde{C}_1(\Omega)\sqrt{\tilde{C}_2(\Omega)\tilde{C}_0}
+2\tilde{C}_2(\Omega)\tilde{C}_0\right)\sqrt{\epsilon_0}
\|\Delta\mf{d}\|_{\mf{L}^2(\Omega)}^2+\frac{1}{4\tilde{C}_3}\|\Delta\mf{d}\|_{\mf{L}^2(\Omega)}^2\right]\\
&\quad  +C(\Omega)(\|\nabla\mf{d}\|_{\mf{L}^2(\Omega)}^2+\|\mf{d}-{\mf{e}}_N\|_{\mathbf{L}^2(\Omega)}^2),
\end{aligned}\end{equation}
where $\tilde{C}_3\geq 2^{-1}$ denotes a constant. Now, choosing $\epsilon_0:=\epsilon_0(\Omega)\in (0,1]$ such that
\begin{equation}\label{choosingjww0219}
\sqrt{\epsilon_0}\tilde{C}_3\left(\tilde{C}_1(\Omega)\sqrt{\tilde{C}_2(\Omega)\tilde{C}_0}
+2\tilde{C}_2(\Omega)\tilde{C}_0\right)\leq \frac{1}{4}, \end{equation}
we get then
\begin{equation}\label{energyinequality}\begin{aligned}&\frac{d}{dt}\mathcal{E}_\delta(t)+\left\|
\mu |\nabla \mf{v}|^2+(\lambda+\mu)|\mm{div}\mf{v}|^2 +\frac{\nu\theta}{2}|\Delta
\mf{d}|^2+\varepsilon\delta\beta\rho^{\beta-2}|\nabla
\rho|^2+{A\varepsilon\gamma}\rho^{\gamma-2}|\nabla\rho|^2\right\|_{L^1(\Omega)}\\
&\leq C(\Omega)(\|\nabla\mf{d}\|_{\mf{L}^2(\Omega)}^2+\|\mf{d}-{\mf{e}}_N\|_{\mathbf{L}^2(\Omega)}^2)\leq C(\Omega)(\mathcal{E}_\delta(t)+1),
\end{aligned}\end{equation}
which, together with Gronwall's inequality, implies that
\begin{equation*}%\label{}
\begin{aligned}
\mathcal{E}_\delta(t)\leq C( {\mathcal{E}}_{\delta,0}, T,\Omega)\quad\mbox{ for any }t\in (0,T).
\end{aligned}\end{equation*}
 Consequently, we can further infer that
\begin{equation}\label{jwwn0221}\begin{aligned}&\mathcal{E}_\delta(t)+\left\|\mu |\nabla \mf{v}|^2
+(\lambda+\mu)|\mm{div}\mf{v}|^2  +\frac{\nu\theta}{2}|\Delta\mf{d}|^2
+\varepsilon\delta\beta\rho^{\beta-2}|\nabla \rho|^2
+{A\varepsilon\gamma}\rho^{\gamma-2}|\nabla\rho|^2\right\|_{L^1(Q_t)} \\
&\leq C( {\mathcal{E}}_{\delta,0},T,\Omega),
\end{aligned}\end{equation}
where $Q_t:=(0,t)\times \Omega$ for any $t\in (0,T)$. Moreover, from \eqref{jww0216}, \eqref{jww0218}
and \eqref{jwwn0221} we get
\begin{equation}\label{jww0222}\begin{aligned}
& \|\nabla \mf{d}\|_{L^4(I,\mathbf{L}^4(\Omega))}^2+ \|\nabla^2 \mf{d}\|_{\mathbf{{L}}^2(\Omega)}+
\|\nabla \mf{v}\|_{L^2(I,\mathbf{L}^2(\Omega))}^2 \leq C( {\mathcal{E}}_{\delta,0}, T,\Omega).
\end{aligned}\end{equation}

Finally, using (\ref{jww0222}), H\"{o}lder's and Sobolev's inequalities, we find from the equation
(\ref{jww0103}) that
\begin{equation}\label{partialtd2}
\begin{aligned}
\|\partial_t\mf{d}\|_{L^{4/3}(I,\mathbf{{L}}^2(\Omega))}\leq &C(\|\mathbf{v}\cdot
 \nabla\mf{d} \|_{L^{4/3}(I,\mathbf{{L}}^2(\Omega))}+\|\Delta\mf{d}
+f_\varepsilon(|\nabla \mathbf{d}|^2)\mf{d}\|_{L^{4/3}(I,\mathbf{{L}}^2(\Omega))}) \\
\leq & C( {\mathcal{E}}_{\delta,0},T, \Omega).
\end{aligned}\end{equation}
Similarly, we can also deduce that
\begin{equation}\label{partialtd}\begin{aligned}
 \|\partial_t\mf{d}\|_{L^{2}(I,(\mathbf{{H}}^{1}(\Omega))^*)}\leq & C( {\mathcal{E}}_{\delta,0}, T, \Omega),
\end{aligned}\end{equation}
where $(\mathbf{H}^{1}(\Omega))^*$ denotes the dual space of $\mathbf{H}^1(\Omega)$.

In addition, when $\Omega=B_R$ with $R\geq 1$, we can show that all
the previous estimates on $(\rho,\mf{v},\mf{d})$ are independent of
$B_R$, except for $\partial_t \mf{d}$. In fact, using
(\ref{jww0216}) and (\ref{jww0218}) for $\Omega=B_1$, and scaling
the spatial variables, we can obtain
\begin{equation*}\label{jwwn0224}
\|\nabla \mf{d}\|_{\mathbf{L}^4(B_R)}^4\leq
C_0(\| \nabla^2 \mf{d}\|_{\mathbf{{L}}^{2}(B_R)}^{2}\|\mf{d}-{\mf{e}}_N\|_{\mathbf{L}^\infty(B_R)}^{2}
+\|\mf{d}-{\mf{e}}_N\|_{\mathbf{L}^4(B_R)}^{4}),  \end{equation*}
and
 \begin{equation*}\begin{aligned}\label{jwwn0225}
 \|\nabla^2 \mf{d}\|_{\mathbf{{L}}^2(B_R)}\leq C_0(\|\Delta
 \mf{d}\|_{\mathbf{{L}}^2(B_R)}+\|\nabla\mf{d}\|_{\mathbf{{L}}^2(B_R)}).
 \end{aligned}\end{equation*}
 Hence, repeating the deduction process of \eqref{energyinequality}, and employing
 the above two inequalities, one can have the following estimate:
 \begin{equation}\label{energyinequalityno}\begin{aligned}
 & \frac{d}{dt}\mathcal{E}_\delta(t)+\Big\|
\mu |\nabla \mf{v}|^2+(\lambda+\mu)|\mm{div}\mf{v}|^2 +\frac{\nu\theta}{2}|\Delta
\mf{d}|^2+\varepsilon\delta\beta\rho^{\beta-2}|\nabla
\rho|^2+{A\varepsilon\gamma}\rho^{\gamma-2}|\nabla\rho|^2\Big\|_{L^1(\Omega)}\\
&\leq C(\|\nabla\mf{d}\|_{\mf{L}^2(\Omega)}^2+\|\mf{d}-{\mf{e}}_N\|_{\mathbf{L}^2(\Omega)}^2),
\end{aligned}\end{equation}
where the constant $C$ is independent of $\Omega=B_R$ for any $R\geq 1$. On the other hand,
using $\eqref{jww0216bounded}$, \eqref{choosingjww0219} and Cauchy-Schwarz's inequality, we
see from (\ref{jww0103}) that
\begin{equation*}\begin{aligned}\label{jwwn0226}
&\frac{1}{2}\frac{d}{dt}\|\mf{d}-\mf{e}_N\|^2_{\mf{L}^2(B_R)}+\theta\|\nabla
\mf{d}\|^2_{\mf{L}^2(B_R)}
=\theta\int_{B_R}\left[f_\varepsilon(|\nabla \mf{d}|^2)(1-d_N)+\frac{1}{2}|
 \mf{d}-\mf{e}_N|^2\mm{div}\mf{v}\right]\mm{d}\mf{x}\\
&\leq \theta\tilde{C}_0\sqrt{\epsilon_0}\|\nabla
\mf{d}\|_{\mf{L}^2(B_R)}^2+ \frac{ \theta\tilde{C}_0\sqrt{\varepsilon_0}}{2}\|
\mm{div}\mf{v}\|_{{L}^2(B_R)}\|\mf{d}-\mf{e}_N\|_{\mf{L}^2(B_R)} \\
&\leq \frac{\theta}{2}\|\nabla \mf{d}\|_{\mf{L}^2(B_R)}^2
+\frac{\mu}{2}\|\mm{div}\mf{v}\|_{{L}^2(B_R)}+C(\theta,\mu^{-1})\|\mf{d}-\mf{e}_N\|_{\mf{L}^2(B_R)}^2.
  \end{aligned}\end{equation*}
Adding the above estimate to \eqref{energyinequality}, we get
\begin{equation*}%\label{}
\begin{aligned}&\left\|
\frac{\mu|\nabla \mf{v}|^2+\theta|\nabla \mf{d}|^2}{2}+(\lambda+\mu)|\mm{div}\mf{v}|^2 +\frac{\nu\theta}{2}|\Delta
\mf{d}|^2+\varepsilon\delta\beta\rho^{\beta-2}|\nabla
\rho|^2+{A\varepsilon\gamma}\rho^{\gamma-2}|\nabla\rho|^2\right\|_{L^1(\Omega)}\\
&+\frac{d}{dt}\left(\mathcal{E}_\delta(t)+\|\mf{d}(t)-\mf{e}_N\|^2_{\mf{L}^2(B_R)}\right)
\leq C\left(\mathcal{E}_\delta(t)+\|\mf{d}(t)-\mf{e}_N\|^2_{\mf{L}^2(B_R)}\right),
\end{aligned}\end{equation*}
which, together with Gronwall's inequality, yields
\begin{equation*}%\label{}
\begin{aligned}&\left\|
\frac{\mu|\nabla \mf{v}|^2+\theta|\nabla \mf{d}|^2}{2}+(\lambda+\mu)|\mm{div}\mf{v}|^2 +\frac{\nu\theta}{2}|\Delta
\mf{d}|^2+\varepsilon\delta\beta\rho^{\beta-2}|\nabla
\rho|^2+{A\varepsilon\gamma}\rho^{\gamma-2}|\nabla\rho|^2\right\|_{L^1(Q_t)}\\
&+\mathcal{E}_\delta(t)+\|\mf{d}(t)-\mf{e}_N\|^2_{\mf{L}^2(B_R)}\leq  C(\mathcal{E}_{\delta,0},T,\|\mf{d}_0-\mf{e}_N\|_{\mf{L}^2(\Omega)}),
\qquad\mbox{for all }\; t\in I.
\end{aligned}\end{equation*}

%To derive the corresponding estimates in the 2D case,
%we shall use the following two special cases of Lemma \ref{Nirenberg}:
%\begin{eqnarray}
%&&\|\nabla {u}\|_{\mf{L}^4(\Omega)}^2\leq
%C(\Omega)(\| \nabla^2 u\|_{\mf{L}^{2}(\Omega)}\|u\|_{L^\infty(\Omega)}+\|u\|_{L^2(\Omega)}^2), \nonumber\\
%&& \label{twqhjf2d}\|{u}\|_{L^4(\Omega)}^2 \leq C(\Omega)(\| \nabla u\|_{\mf{L}^{2}(\Omega)}
%\|u\|_{L^2(\Omega)}+\|u\|_{L^2(\Omega)}^2)\quad \mbox{ for }u\in H^2(\Omega).
%\end{eqnarray}
%Thus, repeating the above deduction procedure, we find that the estimates \eqref{jwwn0221}--\eqref{partialtd}
%still hold in the 2D case. Moreover, the constants in \eqref{jwwn0221} and \eqref{jww0222} for the 2D case
%are also independent of $\Omega=B_R$ for any $R\geq 1$.
%% Of course, the estimate constants will vary, since the dimension of space vary.

Summing up the above estimates, we conclude that
%%%%%%%%%%%%%%%%%%%%%%%%%%%%%%%%%%%%%%%%
\begin{pro}\label{pro:energy}
Let $N=2$ or $3$, and $\Omega\subset \mathbb{R}^N$ be a bounded domain and $\epsilon_0:=\epsilon_0(N,\Omega)\in (0,1]$
satisfy \eqref{choosingjww0219}. Then  the initial-boundary value problem (\ref{jww0101})--(\ref{jwwn0306})
enjoys the following {a priori} estimates, provided the initial data ${d}_{0N}$ satisfies
$1-d_{0N}<\epsilon_0$ and $|\mf{d}_0|\leq 1$.
\begin{eqnarray}
\label{equal1}&& |\mf{d}|\leq 1,\;\mbox{ in particular, }|\mf{d}|=1 \mbox{ if }|\mf{d}_0|=1\mbox{ for }N=2,\\
&& \|\mf{d}-{\mf{e}}_N\|_{\mathbf{L}^\infty(\Omega)}\leq C_0\sqrt{\epsilon_0},\\
&&\label{jww0227}\sup_{t\in I}\|( \mf{d}(t)-\mf{e}_N,\nabla
\mathbf{d}(t))\|_{\mf{L}^2(\Omega)}+\|\nabla^2 \mf{d}\|_{\mf{L}^2(Q_T)}+\|\nabla
\mf{d}\|_{\mf{L}^4(Q_T)}\leq C( {\mathcal{E}}_{\delta,0},T,\Omega),\\
 &&\label{verh1}\sup\limits_{t\in I}(\|\sqrt{\rho}{\bf
v}(t)\|_{\mathbf{L}^2(\Omega)}+\|\rho(t)\|_{{L}^\gamma(\Omega)}) +\|\nabla
{\mathbf{v}}\|_{L^2(I,\mf{L}^2(\Omega))}\leq C( {\mathcal{E}}_{\delta,0},T,\Omega), \\
&&\label{jww0229}\|\partial_t\mathbf{d}\|_{L^{4/3}(I,\mathbf{L}^2(\Omega))}
+\|\partial_t\mathbf{d}\|_{L^{2}(I,(\mf{H}^{1}(\Omega))^*)}\leq C( {\mathcal{E}}_{\delta,0},\Omega,T).
\end{eqnarray}
In particular, if $\Omega=B_R$ with $R\geq 1$, then $\epsilon_0$ can be chosen to be independent
of the domain $\Omega=B_R$ for any $R\geq 1$, and the constant $C( {\mathcal{E}}_{\delta,0},T,\Omega)$
 in \eqref{jww0227} and \eqref{verh1} can be replaced by a constant
 $C({\mathcal{E}}_{\delta,0},T,\|\mf{d}_0-\mf{e}_N\|_{\mf{L}^2(\Omega)})$ independent of $B_R$.
\end{pro}

\section{Strong solvability of sub-systems in the third approximate problem}\label{subsystem}

Before proving the unique solvability of the third approximate problem, we shall introduce two preliminary results.
The first result is concerned with the global  solvability of the Neumann problem on the equation (\ref{jww0101})
for given $\mf{v}$.
%%%%%%%%%%%%%%%%%%%%%%%%%%%%%%%%%%%%%%%%%%
\begin{pro}\label{pro:0401}
Let $0<\alpha<1$, $\Omega$ be a bounded domain of class $C^{2,\alpha}$,
and $\rho_0:=\rho_0(\mathbf{x})$ satisfy
 \begin{equation*}\rho_0\in W^{1,\infty}(\Omega),
 \qquad 0<\underline{\rho}\leq\rho_0\leq \bar{\rho}<\infty. \end{equation*}
 Then, there exists a unique mapping
\begin{equation*}\begin{aligned}
\mathscr{S}_{\rho_0}:
L^{\infty}(I,\mathbf{W}^{1,\infty}_0(\Omega))\rightarrow\ C^0(\bar{I},H^{1}(\Omega)),
\end{aligned}\end{equation*}
such that
\begin{enumerate}
  \item[ \qquad (1)] $\mathscr{S}_{\rho_0}(\mathbf{v})$ belongs to the function class
  \begin{equation*}\label{n0404}
\mathcal{R}_T:=\left\{\rho\,|\;\rho\in L^2(I,W^{2,q}(\Omega))\cap C^0(\bar{I},W^{1,q}(\Omega)),\;
\partial_t\rho\in L^2(I,L^q(\Omega))\right\},\;\; 1<q<\infty.
\end{equation*}
  \item[ \qquad (2)]  The function $\rho=\mathscr{S}_{\rho_0}(\mathbf{v})$
  satisfies the following initial-boundary problem:
  \begin{equation}\label{jww0309}   \left\{\begin{array}{ll}
  \partial_t\rho+\mathrm{div}(\rho\mathbf{v})=\varepsilon \Delta \rho &\mbox{ a.e. in } Q_T,\\[0.4em]
 \rho(\mathbf{x},0)=\rho_0  &\mbox{ a.e. in }\Omega,\\[0.4em]
 \nabla\rho\cdot\mf{n}|_{\partial\Omega}=0,&\mbox{ in the sense of traces a.e. in }I.
     \end{array}\right.  \end{equation}
\item[ \qquad (3)] $\mathscr{S}_{\rho_0}(\mathbf{v})$ is pointwise bounded, i.e.,
 \begin{eqnarray}\label{0404}
 \underline{\rho}e^{-\int_0^t\|\mf{v}\|_{\mf{W}^{1,\infty}(\Omega)}\mm{d}\tau}
 \leq \mathscr{S}_{\rho_0}(\mf{v})(\mf{x},t)\leq \bar{\rho}
 e^{\int_0^t\|\mf{v}(\tau)\|_{\mf{W}^{1,\infty}(\Omega)}\mm{d}\tau},\ t\in \bar{I},
 \quad\mbox{for a.e. }\mf{x}\in \Omega.\end{eqnarray}
 \item[ \qquad (4)] If $\|\mathbf{v}\|_{L^\infty(I,\mathbf{W}^{1,\infty}(\Omega))}\leq\kappa_v$, then
\begin{eqnarray}
&&\|\mathscr{S}_{\rho_0}(\mathbf{v})\|_{L^\infty(I_t,H^{1}(\Omega))}\leq C(\Omega)\|\rho_0
\|_{H^1(\Omega)}e^{\frac{C(\Omega)}{2\varepsilon}(\kappa_{v}+\kappa^2_{v})t},  \nonumber \\
&& \|\nabla^2\mathscr{S}_{\rho_0}(\mathbf{v})\|_{\mf{L}^2(Q_t)}\leq
\frac{C(\Omega)}{\varepsilon}\sqrt{t} \|\rho_0\|_{H^1(\Omega)}\kappa_{v}
e^{\frac{C(\Omega)}{2\varepsilon}(\kappa_{v}+\kappa^2_{v})t}, \label{0405} \\
&& \|\partial_t\mathscr{S}_{\rho_0}(\mathbf{v})\|_{L^2(Q_t)}\leq
C(\Omega)\sqrt{t}\|\rho_0\|_{H^1(\Omega)} \kappa_{v}
e^{\frac{C_1}{2\varepsilon}(\kappa_{v}+\kappa^2_{v})t} \nonumber\end{eqnarray}
for any $t\in \bar{I}$, where $I_t:=(0,t)$, and $Q_t:=\Omega\times I_t$.
\item[ \qquad (5)]  $\mathscr{S}_{\rho_0}(\mathbf{v})$ depends continuously  on $\mf{v}$,
  i.e.,
\begin{equation*}\begin{aligned}\label{0409}\|[\mathscr{S}_{\rho_0}
(\mathbf{v}_1)-\mathscr{S}_{\rho_0}(\mathbf{v}_2)](t)\|_{L^2(\Omega)}\leq
G(\kappa_v,\varepsilon,T)t\|\rho_0\|_{H^1(\Omega)}
\|\mathbf{v}_1-\mathbf{v}_2\|_{L^\infty(I_t,\mathbf{W}^{1,\infty}(\Omega))},\end{aligned}
\end{equation*}
\begin{equation*}\begin{aligned}\label{n04091}\|\partial_t[\mathscr{S}_{\rho_0}
(\mathbf{v}_1)-\mathscr{S}_{\rho_0}(\mathbf{v}_2)](t)\|_{L^2(Q_t)}\leq
G(\kappa_v,\varepsilon,T)\sqrt{t}\|\rho_0\|_{H^1(\Omega)}
\|\mathbf{v}_1-\mathbf{v}_2\|_{L^\infty(I_t,\mathbf{W}^{1,\infty}(\Omega))}\end{aligned}
\end{equation*}
for any $t\in \bar{I}$, and for any
$\|\mathbf{v}_1\|_{L^\infty(I,\mathbf{W}^{1,\infty}(\Omega))}\leq
\kappa_v$ and $\|\mathbf{v}_2\|_{L^\infty(I,\mathbf{W}^{1,\infty}(\Omega))}\leq
\kappa_v$. The constant $G$ is nondecreasing in the first
variable and may depend on $\Omega$.
\end{enumerate}
\end{pro}

\begin{pf} Please refer to \cite[Proposition 7.39]{NASII04} and \cite[Proposition 3.1]{JFJSWDH}\hfill $\Box$
\end{pf}

The second result is on the local solvability of the Neumann problem for the system (\ref{jww0103})
with given $\mf{v}$.
%%%%%%%%%%%%%%%%%%%%%%%%%%%%%%%%%%%
\begin{pro}\label{pro:0402}
Let
$0<\alpha<1$, $\Omega\in \mathbb{R}^N$ be a bounded domain of class
$C^{2,\alpha}$, $\mf{d}_0:=\mf{d}_0(\mf{x})\in \mf{H}^2(\Omega)$, $K>0$ and
$$\mf{v}\in \mathbb{V}_K:=\left\{\mf{v}~\bigg|~\|\mf{v}\|_{\mathbb{V}}
:=\left(\|\mathbf{v}\|_{C^0 (\bar{I},\mathbf{H}^{2}(\Omega))}^2
+\|\partial_t\mf{v}\|_{L^{2}(I,\mathbf{H}^{1}(\Omega))}^2\right)^\frac12\leq K\right\}.$$
Then there exist a finite time
\begin{equation*}
{T_d^K}:=h_1(K,\|\mf{d}_0\|_{\mathbf{H}^2(\Omega)},\Omega)\in (0, \min\{1,T\}),\end{equation*}
and a corresponding unique mapping
\begin{equation*}
\mathscr{D}^K_{\mathbf{d}_0}:
\mathbb{V}_K (\mbox{with }T^K_d\mbox{ in place of }T)\to C^0(\bar{I}^K_d,\mf{H}^2(\Omega)),
\end{equation*}
where $h_1$ is nonincreasing in its first two variables and
${Q}_{T^K_d}:={\Omega}\times I^K_d:=\Omega\times (0,T^K_d)$, such that
\begin{enumerate}
  \item[\qquad (1)] $\mathscr{D}^K_{\mathbf{d}_0}(\mathbf{v})$ belongs to the following function class
\begin{equation}\label{0465}\begin{aligned}
\mathcal{R}_{T_d^K}:=\{\mf{d}~|~ & \mf{d}\in L^2(I_d^K,\mf{H}^3(\Omega)),\ \partial_{tt}^2\mf{d}
\in {L^2(I_d^K,(\mf{H}^{1}(\Omega))^*)},\\
& \partial_t\mf{d}\in {C^0(\bar{I}_d^K,\mf{L}^2(\Omega))}\cap
{L^2(I_d^K,\mf{H}^1(\Omega))} \}.\end{aligned}\end{equation}
  \item[\qquad (2)]  $\mf{d}=\mathscr{D}^K_{\mathbf{d}_0}(\mathbf{v})$ satisfies the following initial-boundary problem:
\begin{equation}\label{jww0317}\left\{\begin{array}{ll}
\partial_t\mathbf{d}+\mathbf{v}\cdot \nabla \mathbf{d}= \theta(\Delta
\mathbf{d}+f_\varepsilon(|\nabla \mathbf{d}|^2)\mf{d})&\mbox{ a.e. in }Q_{T_d^K},  \\[0.4em]
\mathbf{d}(\mf{x},0)=\mathbf{d}_0&\mbox{ a.e. in }\Omega,\\[0.4em]
(\mf{n}\cdot\nabla \mathbf{d})|_{\partial\Omega}=\mathbf{0}&\mbox{ in the sense of trace a.e. in }I^K_d,
\end{array}\right.\end{equation} where $f_\varepsilon\geq 0$ is defined by (\ref{appfun}).
   \item[\qquad (3)]  $\mathscr{D}^K_{\mathbf{d}_0}(\mathbf{v})$ enjoys the following estimate:
\begin{equation*}\label{0455}
\|\mathscr{D}^K_{\mathbf{d}_0}(\mathbf{v})\|_{\mathbb{D}^*}\leq
C(K,\|\mf{d}_0\|_{\mathbf{H}^2(\Omega)},\Omega),
\end{equation*}
where
$$\|\cdot\|_{\mathbb{D}^*}:=\left(\|\partial_t\cdot\|_{L^2(I_d^K,\mf{H}^1(\Omega))}^2+\|
\partial_t\cdot\|_{L^\infty(I_d^K,\mf{L}^2(\Omega))}^2+\|\cdot\|_{L^\infty(I_d^K,\mf{H}^2(\Omega))}^2
+\|\cdot\|_{L^2(I_d^K,\mf{H}^3(\Omega))}^2\right)^\frac12.$$
Moreover (recalling (\ref{equal1}) and (\ref{jww0205})),
\begin{eqnarray}
&&|\mf{d}|\leq 1,\;\mbox{ in particular, } |\mf{d}|\equiv1 \mbox{ if }|\mf{d}_0|\equiv 1\mbox{ for }N=2;\nonumber\\
&&  d_{i}(t)\geq \underline{d}_{0i}\;\;\mbox{ for any }t,\;\;\mbox{ if }\; {d}_{0i}
\geq \underline{d}_{0i}\geq 0\ \mbox{ where }1\leq i\leq N.\nonumber
\end{eqnarray}
  \item[\qquad (4)]   $\mathscr{D}^K_{\mathbf{d}_0}(\mathbf{v})$ continuously depends on
  $\mf{v}$ in the following sense:
  \begin{equation*}\label{n0456}
  \begin{aligned}&
\|[\mathscr{D}^K_{\mathbf{d}_0}(\mathbf{v}_1)-\mathscr{D}^K_{\mathbf{d}_0}
(\mathbf{v}_2)](t)\|_{\mathbf{L}^2(\Omega)} +
\|\nabla(\mathscr{D}^K_{\mathbf{d}_0}(\mathbf{v}_1)-\mathscr{D}^K_{\mathbf{d}_0}
(\mathbf{v}_2))\|_{\mf{L}^2(Q_{t})}  \\
& \quad \leq \sqrt{t}C(K,\|\mf{d}_0\|_{\mathbf{H}^2(\Omega)},\Omega)
\|\mathbf{v}_1-\mathbf{v}_2\|_{L^\infty(I_t,\mf{L}^2(\Omega))}
\qquad\mbox{for any }\; t\in \bar{I}_d^K.
\end{aligned}\end{equation*}
\end{enumerate}
\end{pro}

\begin{pf} The above results on the two-dimensional non-homogeneous Derichlet problem
have been shown under the higher regularity condition $\mf{d}_0\in\mathbf{H}^3(\Omega)$ in \cite{JFJSWDH},
where the appeared constants depend on $\|\nabla^3\mf{d}_0\|_{\mf{L}^2(\Omega)}$.
For the multi-dimensional Neumann problem considered here, we can obtain the same results
by slightly modifying the proof in \cite[Proposition 3.2]{JFJSWDH} as follows.

 (1) We utilize a semi-discrete Galerkin method from \cite[Proposition 7.39]{NASII04}
 to construct a solution of the linearized Neumann problem,
 and a semi-discrete Galerkin method from \cite{ELGP}
 to construct a solution of the linearized non-homogeneous Derichlet problem in \cite{JFJSWDH}.

(2) We make use of the elliptic estimate (\ref{jww0218}) for the Neumann boundary condition to replace
   the corresponding elliptic estimate for the non-homogeneous boundary condition in \cite{JFJSWDH}.
 This is the reason why the constant $C$ in Proposition \ref{pro:0402} is independent of
 $\|\nabla^3 \mf{d}_0\|_{\mf{L}^2(\Omega)}$.
 %%% It should be remarked here that the term on gradient in the elliptic
%%% estimate (\ref{jww0218}) will not result in any big troubles.

(3) We use the following three-dimensional interpolation inequalities in Lemma \ref{Nirenberg}
\begin{eqnarray}\label{jwwn032012}
&& \| {u}\|_{L^\infty(\Omega)}\leq
C(\Omega)(\|\nabla^2 u\|_{{L}^{2}(\Omega)}^{\frac{3}{4}}\|u\|_{L^2(\Omega)}^{\frac{1}{4}}+\|u\|_{L^2(\Omega)})
\quad\mbox{for }u\in H^2(\Omega), \\
&& \label{jww0208jww}  \|{u}\|_{L^4(\Omega)}\leq
C(\Omega)(\| \nabla u\|_{{L}^{2}(\Omega)}^{\frac{3}{4}}\|u\|_{L^2(\Omega)}^{\frac{1}{4}}+\|u\|_{L^2(\Omega)})
\quad\mbox{for }u\in H^1(\Omega), \\
&& \label{jww0211}   \|{u}\|_{L^4(\Omega)}\leq \sqrt{2}(\|\nabla u\|_{{L}^{2}(\Omega)}^{\frac{3}{4}}
\|u\|_{L^2(\Omega)}^{\frac{1}{4}}\quad\mbox{for }u\in H^1_0(\Omega)\quad\mbox{(cf. \cite[Lemma 3.5]{NSRT1})},
\end{eqnarray}
to replace the following two-dimensional interpolation inequalities
  \begin{eqnarray}\label{jwwn0320}
&&\| {u}\|_{L^\infty(\Omega)}^2\leq
C(\Omega)\| u\|_{{H}^{2}(\Omega)}\|u\|_{L^2(\Omega)}\quad \mbox{ for }u\in H^2(\Omega),\\
&&\|{u}\|_{L^4(\Omega)}^2\leq
C(\Omega)\|  u\|_{{H}^{1}(\Omega)}\|u\|_{L^2(\Omega)}\quad\mbox{ for }u\in H^1(\Omega),\\
&&\label{jwwn0322}\|{u}\|_{L^4(\Omega)}^4\leq \sqrt{2}\| \nabla u\|_{{L}^{2}(\Omega)}^{2}
\|u\|_{L^2(\Omega)}^{2}\quad \mbox{ for }u\in H^1_0(\Omega).
\end{eqnarray}
The purpose of using (\ref{jwwn0320})--(\ref{jwwn0322}) above in \cite{JFJSWDH} is to construct the term $t^\alpha$
with $\alpha\in (0,1)$. When we use \eqref{jwwn032012}--(\ref{jww0211}) in showing solvability of
the three-dimensional problem, we can still construct the term $t^\alpha$ with possible different
$\alpha\in (0,1)$, with the help of Young's inequality in addition.

Finally, we remark that the conclusions in Proposition \ref{pro:0402} still hold in the 3D case, provided
$f_\varepsilon( {x})\equiv  {x}$.    \hfill $\Box$
\end{pf}

\section{Unique solvability of the third approximate problem}\label{sec:jww03}
In order to obtain a weak solution of the initial-boundary problem (\ref{0101})--(\ref{0104}),
we first show the existence of solutions to the third
approximate problem of the original three-dimensional problem (\ref{0101})--(\ref{0104}):
\begin{eqnarray} && \label{n301}
\partial_t\rho+\mathrm{div}(\rho\mathbf{v}) =\varepsilon\Delta\rho ,   \\
&& \label{n302}\partial_t\mathbf{d}+\mathbf{v}\cdot \nabla \mathbf{d}
= \theta(\Delta \mathbf{d}+f_\varepsilon(|\nabla \mathbf{d}|^2)\mf{d}),  \\[1mm]
&& \begin{aligned}\label{n303}&\int_\Omega
(\rho\mathbf{v})(t)\cdot\mathbf{\mathbf{\Psi}}\mathrm{d}\mathbf{x}
-\int_{\Omega}\mathbf{m}_0\cdot\mathbf{\mathbf{\Psi}}\mathrm{d}\mathbf{x} \\
& = \int_0^t\int_\Omega\bigg[\mu\Delta\mathbf{v}+(\mu+\lambda)\nabla
\mm{div}\mf{v}-A\nabla \rho^\gamma -\delta \nabla
\rho^\beta-\varepsilon(\nabla\rho\cdot\nabla \mathbf{v})  \\
&\qquad\qquad -\mathrm{div}(\rho\mathbf{v}\otimes\mathbf{v})
-\nu\mathrm{div}\left(\nabla \mathbf{d}\otimes\nabla\mathbf{d} -\frac{|\nabla\mf{d}|^2
\mathbb{I}}{2}\right)\bigg]\cdot\mathbf{\Psi}\mathrm{d}\mathbf{x}\mathrm{d}s
 \end{aligned} \end{eqnarray}
 for all $t\in  I$  and any $\mathbf{\Psi}\in \mathbf{X}_n$, with boundary conditions
\begin{eqnarray}\nabla\rho\cdot\mf{n}|_{\partial\Omega}=0,\quad
\mathbf{v}|_{\partial \Omega}=\mathbf{0},\quad
(\mf{n}\cdot\nabla\mathbf{d})|_{\partial\Omega}=\mathbf{0},\end{eqnarray} and
modified initial data
\begin{eqnarray}\label{n0305}
&&\rho(\mathbf{x},0)=\rho_0\in W^{1,\infty}(\Omega),
\quad 0<\underline{\rho}\leq\rho_0\leq \bar{\rho}<\infty,  \\
&&\label{n0306}\mathbf{d}(\mathbf{x},0)=\mathbf{d}_0\in\mf{H}^{2}({\Omega}),\quad
 \mathbf{v}(\mathbf{x},0)=\mf{v}_0\in \mathbf{X}_n,\end{eqnarray}
where  $\mf{n}$ denotes the outward normal to $\partial \Omega$,
and $\varepsilon,\delta, \beta, \underline\rho, \bar{\rho}>0$ are constants.

 Here we briefly introduce the finite dimensional space $\mf{X}_n$.
We know from \cite[Section 7.4.3]{NASII04} that there exist
countable sets
\begin{equation*}\begin{aligned}&\{\lambda_i\}_{i=1}^\infty,\
0<\lambda_1\leq \lambda_2\leq \cdots ,\mbox{ and }  \\
&\{\mathbf{\Psi}_i\}_{i=1}^\infty\subset\mathbf{W}_0^{1,p}(\Omega)\cap\mathbf{W}^{2,p}(\Omega),\
1\leq p<\infty,
\end{aligned}\end{equation*}
such that
\begin{equation*} \begin{aligned}
&-\mu\Delta\mathbf{\Psi}_i-(\mu+\lambda)\nabla
\mm{div}\mathbf{\Psi}_i=\lambda_i\mathbf{\Psi}_i,\qquad i=1,2,\cdots ,
\end{aligned}\end{equation*}
and $\{\mathbf{\Psi}_i\}_{i=1}^\infty$ is an orthonormal basis in
$\mf{L}^2(\Omega)$ and an orthogonal basis in $\mf{H}_0^1(\Omega)$
with respect to the scalar product
$\int_\Omega[\mu\partial_j\mf{u}\cdot\partial_j\mf{v}
+(\mu+\lambda)\mm{div}\mf{u}\,\mm{div}\mf{v}]\mm{d}\mf{x}$.
We define a $n$-dimensional Euclidean space $\mf{X}_n$ with scalar
product $<\cdot,\cdot >$ by
\begin{equation*}\begin{aligned}
\mf{X}_i=\mm{span}\{\mathbf{\Psi}_i\}_{i=1}^n,\quad
<\mf{u},\mf{v}> =\int_\Omega\mf{u}\cdot\mf{v}\mm{d}\mf{x},\quad \mf{u},\mf{v}\in\mf{X}_n,
\end{aligned}\end{equation*}
and denote by $\mathscr{P}_n$ the orthogonal projection of $\mf{L}^2(\Omega)$ onto $\mf{X}_n$.
\subsection{Local existence}
With the help of Proposition \ref{pro:0401} and \ref{pro:0402},
one can  establish the local existence of a unique solution to the third approximate
problem (\ref{n301})--(\ref{n0306}) by a contraction mapping argument.
To this purpose, we rewrite the approximate momentum equations (\ref{n303}) as an operator form.

Given
\begin{equation*}\begin{aligned}\label{0504}
\rho\in C^0(\bar{I},L^1(\Omega)),
\quad \partial_t \rho\in L^1(\Omega_{T}),\quad  \underset{(\mathbf{x},t)\in \Omega_T}
{\mathrm{ess\ inf}}\,\rho (\mathbf{x},t)\geq \underline{\rho}>0, \end{aligned}\end{equation*}
we define, for all $t\in \bar{I}$, that
\begin{equation*}\begin{aligned}
\mathscr{M}_{\rho(t)}:\; \mathbf{X}_n \to \mathbf{X}_n
\end{aligned}\end{equation*}
by
$$\ <\mathscr{M}_{\rho(t)}\mathbf{v},\,
\mathbf{w}>:=\int_{\Omega}\rho(t)\mathbf{v}\cdot\mathbf{w}\mathrm{d}x,\quad
\mathbf{v},\mathbf{w} \in \mathbf{X}_n.$$
Recall that all norms on $\mathbf{X}_n$ are equivalent, in particular,
\begin{equation}\label{0506}
\mathbf{W}^{k_1,p_1}(\Omega)\mbox{ and }
(\mathbf{W}^{k_2,p_2}_0(\Omega))^*\mbox{-norms are equivalent on }\mathbf{X}_n,\end{equation}
where $(\mathbf{W}^{k_2,p_2}_0(\Omega))^*$ denotes the dual space of
$(\mathbf{W}^{k_2,p_2}_0(\Omega))$, $k_1$ and $k_2$ are integers, and
$0\leq k_2<\infty$, $1\leq p_2\leq \infty$, $0\leq k_1\leq1$
and $1\leq p_1\leq \infty$ (or $k_1=2$, $1\leq p_1<\infty$).
Note that this property of equivalent norms plays an important role in the estimates of velocity $\mf{v}$.

It is easy to see that
\begin{equation}\begin{aligned}\label{0507}
\|\mathscr{M}_{\rho(t)}\|_{\mathscr{L}(\mathbf{X}_n,\mathbf{X}_n)}\leq c(n)\int_\Omega\rho(t)
\mathrm{d}x,\quad t\in \bar{I}.\end{aligned}\end{equation}
On the other hand, we easily verify that $\mathscr{M}^{-1}_{\rho(t)}$ exists for all
$t\in \bar{I}$ and
\begin{equation}\label{0508}
\|\mathscr{M}^{-1}_{\rho(t)}\|_{\mathscr{L}(\mathbf{X}_n,\mathbf{X}_n)}
\leq \frac{1}{\underline{\rho}},\end{equation}
where $\mathscr{L}(\mathbf{X}_n,\mathbf{X}_n)$ denotes the set of
all continuous linear operators mapping $\mathbf{X}_n$ to $\mathbf{X}_n$.
 By virtue of (\ref{0507}) and (\ref{0508}), we have
\begin{equation}\begin{aligned}\label{0509}\|\mathscr{M}^{-1}_{\rho(t)}\mathscr{M}_{\rho_1(t)}
\mathscr{M}^{-1}_{\rho(t)} \|_{\mathscr{L}(\mathbf{X}_n,\mathbf{X}_n)}\leq\frac{c(n)}
{\underline{\rho}^2}\|\rho_1(t)\|_{L^1(\Omega)},\quad t\in
\bar{I}.\end{aligned}\end{equation}

 Next, we shall look for $T_n^*\subset (0,T^{{K}}_d]$ and
\begin{equation*}
\mathbf{v}\in \mathbb{A}:=\{\mf{v}\in C(\bar{I}_n^*,
\mathbf{X}_n)~|~\partial_t\mf{v}\in L^2(I_n^*,\mf{X}_n)\}, \quad
I_n^*:=(0,T_n^*)\subset (0,T^{{K}}_d)\end{equation*}
with
$\|\mathbf{v}\|_{C(\bar{I}_n^*,\mathbf{H}^2(\Omega))}
+\|\partial_t\mathbf{v}\|_{L^2({I}_n^*, \mathbf{H}^1(\Omega))}\leq {K}$
for some ${K}$, such that   %%%% satisfying
\begin{equation}\begin{aligned}\label{0513}
&\int_\Omega(\rho\mathbf{v})(t) \cdot\mathbf{\mathbf{\Psi}}\mathrm{d}\mathbf{x}
-\int_{\Omega}\mathbf{m}_0\cdot\mathbf{\mathbf{\Psi}}\mathrm{d}\mathbf{x} \\
&= \int_0^t\int_\Omega\bigg[\mu\Delta\mathbf{v}+(\mu+\lambda)\nabla
\mm{div}\mf{v}-A\nabla \rho^\gamma -\delta \nabla
\rho^\beta-\varepsilon(\nabla\rho\cdot\nabla \mathbf{v}) \\
& \qquad\qquad -\mathrm{div}(\rho\mathbf{v}\otimes\mathbf{v})
-\nu\mathrm{div}\left(\nabla \mathbf{d}\otimes\nabla\mathbf{d}-\frac{|\nabla
\mf{d}|^2\mathbb{I}}{2}\right)\bigg]\cdot\mathbf{\Psi}\mathrm{d}\mathbf{x}\mathrm{d}s
\end{aligned}\end{equation}
for all $t\in [0,T_n]$ and any $\mathbf{\Psi}\in\mathbf{X}_n$,
where $\rho(t)=[\mathscr{S}_{\rho_0}(\mathbf{v})](t)$ is the solution of the
problem (\ref{jww0309}) constructed in Proposition \ref{pro:0401},
$\mathbf{d}(t)=\mathscr{D}_{\mf{d}_0}^{{K}}(\mathbf{v})(t)$
is the solution of the problem (\ref{jww0317}) constructed in
Proposition \ref{pro:0402}. By the regularity of
$(\rho,\mf{d})$ in Propositions \ref{pro:0401} and \ref{pro:0402}, and
the operator $\mathscr{M}_{\rho(t)}$, the equations (\ref{0513}) can be rephrased as
\begin{equation*}\label{n0512}
 {\mathbf v}(t) =\mathscr{M}^{-1}_{[\mathscr{S}_{\rho_0}({\mf{
v}})](t)}\bigg((\mathscr{P}\mathbf{m}_0+\int^t_0
\mathscr{P}[{\mathscr{N}(\mathscr{S}_{\rho_0}({\bf v})},{\mf{ v}},
\mathscr{D}_{\mf{d}_0}^{{K}}(\mathbf{v}))]\mathrm{d}s\bigg)
\end{equation*}
with $\mf{m}_0=(\rho\mf{v})(0)$, where $\mathscr{P}:=\mathscr{P}_n$
is the orthogonal projection of $\mathbf{L}^2(\Omega)$ to $\mathbf{X}_n$, and
\begin{equation*}\label{0515}\begin{aligned}
&\mathscr{N}(\rho,\mathbf{v},\mathbf{d})=\mu\Delta\mathbf{v}+(\mu+\lambda)\nabla
\mm{div}\mf{v}-A\nabla \rho^\gamma -\delta \nabla
\rho^\beta-\varepsilon(\nabla\rho\cdot\nabla \mathbf{v})\\
&\quad\quad\quad\quad\quad -\mathrm{div}(\rho\mathbf{v}\otimes\mathbf{v})
-\nu\mathrm{div}\left(\nabla \mathbf{d}\otimes\nabla
\mathbf{d}-\frac{|\nabla \mf{d}|^2\mathbb{I}}{2}\right).
\end{aligned}
\end{equation*}
Moreover, one has
\begin{equation}\label{n0413nn}\begin{aligned}
\partial_t\mathbf{v}(t) & =\mathscr{M}_{[\mathscr{S}_{\rho_0}({\mf{v}})](t)}^{-1}
\mathscr{M}_{\partial_t[\mathscr{S}_{\rho_0}({\mf{v}})](t)}
\mathscr{M}_{[\mathscr{S}_{\rho_0}({\mf{v}})](t)}^{-1}
\bigg\{\mathscr{P}\mathbf{m}_0 +\int_0^t[\mathscr{P}\mathscr{N}(\mathscr{S}_{\rho_0}(\mathbf{v}),
\mathbf{v},\mathscr{D}_{\mf{d}_0}(\mf{v}))](s)\mathrm{d}s\bigg\} \\
&\quad +\mathscr{M}_{[\mathscr{S}_{\rho_0}({\mf{ v}})](t)}^{-1}[\mathscr{P}
\mathscr{N}(\mathscr{S}_{\rho_0}(\mathbf{v}),\mathbf{v},\mathscr{D}^{{K}}_{\mf{d}_0}(\mf{v}))](t).
\end{aligned}\end{equation}

The authors in \cite[Section 4.3]{JFJSWDH} have established the
local existence of the problem (\ref{n301})--(\ref{n0306}) with
non-homogenous boundary condition in place of the Neumann boundary
condition for the system (\ref{n302}) by using a contraction mapping argument.
In view of the proof in \cite{JFJSWDH} and the previous preliminary results,
we can immediately find that the results in \cite[Section 4.3]{JFJSWDH}
 can be directly generalized to our problem (\ref{n301})--(\ref{n0306})
without essential changes in arguments. Thus, we have the following conclusion.
%%%%%%%%%%%%%%%%%%%%%%%%%%%%%%%%%%%%%%%%%%
\begin{pro}\label{localtime}
There exist $K:=K(\|\mathscr{P}\mathbf{m}_0\|_{\mathbf{X}_n},\|\mathbf{v}(0)\|_{\mathbf{X}_n},
\underline{\rho}^{-1})>0$ and $T_n^*\leq T$, such that
\begin{eqnarray*}\label{n0530}
\displaystyle\mathcal{T}:\ \mathbb{A}\rightarrow \mathbb{A},\ \mathcal{T}(\mathbf{w})
:=\mathscr{M}^{-1}_{[\mathscr{S}(\mathbf{w})]}
\left\{\mathscr{P}\mathbf{m}_0+\displaystyle\int_0^t[\mathscr{P}\mathscr{N}(\mathscr{S}(\mathbf{w}),
\mathbf{w},\mathscr{D}(\mathbf{w}))](s)\mathrm{d}s\right\}, \end{eqnarray*}
maps
\begin{equation*}B_{K,\tau_0}=\big\{\mathbf{w}\in \mathbb{A}~|~
 \|\mathbf{w}\|_{C(\bar{I}_{\tau_0},\mathbf{X}_n)}
+\|\partial_t\mathbf{w}\|_{L^2({I}_{\tau_0},\mathbf{X}_n)}\leq
K\big\}\subset C(\bar{I}_{{T}_n^*},\mathbf{X}_n)\end{equation*}
into itself and is contractive for any $0<\tau_0\leq T_n^*$,
 where one can take $\mf{X}_n=\mf{H}^1_0(\Omega)$, $T_n^*$ has the form
\begin{equation}\label{n0540}
0<{T}_n^*={h}_2({\bar{\rho},\|\rho_0\|_{H^1(\Omega)},\|\mathbf{d}_0\|_{\mf{H}^2(\Omega)},
\underline{\rho},K,T,n}),
\end{equation}
and $h_2$ is nonincreasing in its first three variables and nondecreasing in the fourth variable.
 \end{pro}

  Therefore, the map $\mathcal{T}$ possesses in $B_{K,{T}_n^*}$ a unique
fixed point $\mathbf{v}$ which satisfies (\ref{0513}). Thus, we have a solution
$(\rho=\mathscr{S}(\mathbf{v}),\mathbf{v},\mathscr{D}(\mathbf{v}))$
which is defined in $ Q_{{T}_n^*}$ and satisfies the
initial-boundary value problem (\ref{n301})--(\ref{n0306}) for each given
$n$. This means that we can find a unique maximal solution
$(\rho_n,\mf{v}_n,\mf{d}_n)$ defined in $[0,T_n)\times \Omega$ for
each given $n$, where $T_n\leq T$.

\subsection{Global existence}\label{sec:0404}

In order to show the maximal time $T_n=T$ for any $n$, it
suffices to derive uniform bounds for $\rho_n$, $\mathbf{v}_n$,
$\mathbf{d}_n$ and $\mathscr{P}_n(\rho_n\mathbf{v}_n)$. However, we
need to impose an additional smallness condition as in (\ref{jww0212}) to get the uniform boundedness of
$\|\mf{d}_n\|_{L^\infty(\bar{I}_n,\mathbf{H}^2(\Omega))}$.
For simplicity, we denote
$$(\rho,\mathbf{v},\mathbf{d},\mathscr{P}\mathbf{m}):=(\rho_n,
\mathbf{v}_n, \mathbf{d}_n,\mathscr{P}_n(\rho_n\mathbf{v}_n)).$$
We mention that in the estimates that follow, the letters $G_1(\ldots)$, $G_2(\ldots)$ and $G(\ldots)$ will
denote various positive constants depending on its variables. %%% and may depend on the physical parameters.

First, one has the energy estimates as in Proposition \ref{pro:energy}. In fact, by virtue of the
regularity of $(\rho,\mf{v},\mf{d})$, we can deduce that $(\rho,\mf{v},\mf{d})$ satisfies
(\ref{0542}) for $N\geq 2$, and (\ref{05420}) and \eqref{05421} for $N=2$.
Then, letting $\epsilon_0$ satisfy \eqref{choosingjww0219} and the initial data $\mf{d}_0$ satisfy
\begin{equation*}\begin{aligned}1-{\epsilon_0(\Omega)} \leq  d_{0N}\leq 1\mbox{ and } |\mf{d}_0|\leq 1,
\end{aligned}\end{equation*}
arguing in the same manner as in the derivation of (\ref{equal1})--(\ref{jww0229}), we get
\begin{eqnarray}&& \nonumber|\mf{d}|\leq 1,\;\mbox{ in particular, }\;
|\mf{d}|=1 \mbox{ if }|\mf{d}_0|=1\;\mbox{ for }N=2,\\
&&\label{123654}\|\mf{d}-\mf{e}_N\|_{\mf{L}^\infty(\Omega)}\leq C_0\sqrt{\epsilon_0},\\
\label{n0544}&&\|\sqrt{\rho}\mathbf{v}\|_{L^\infty(I_n,\mathbf{L}^{2}(\Omega))}\leq
G({\mathcal{E}}_{\delta,0},T,\Omega),  \\
&& \label{n0544n}   \|\mathbf{v}\|_{L^2(I_n,\mathbf{H}^{1}(\Omega))}\leq
G({\mathcal{E}}_{\delta,0},T,\Omega),\\&&\label{n054412}   \|\nabla^2 \mf{d}\|_{\mf{L}^2(Q_{T_n})}+\|\nabla
\mf{d}\|_{\mf{L}^4(Q_{T_n})}+
 \|\nabla\mathbf{d}\|_{L^\infty(I_n,\mathbf{L}^{2}(\Omega))}\leq
G({\mathcal{E}}_{\delta,0}, T,\Omega),  \\
&&\|\rho\|_{L^\infty(I_n,L^\gamma(\Omega)}+\|\partial_t\mathbf{d}\|_{L^{4/3}(I_n,\mathbf{L}^2(\Omega))}
+\|\partial_t\mathbf{d}\|_{L^{2}(I_n,(\mf{H}^{1}(\Omega))^*)}\leq
G({\mathcal{E}}_{\delta,0}, T,\Omega).\nonumber
\end{eqnarray}
In particular, if $\Omega=B_R$ with $R\geq 1$, the constant $\epsilon_0$ can be chosen to be independent
of $\Omega$, and the above constant $G({\mathcal{E}}_{\delta,0}, T,\Omega)$ can be replaced by
a constant $G({\mathcal{E}}_{\delta,0},T, \|\mf{d}_0-\mf{e}_N\|_{\mf{L}^2(\Omega)})$ independent of $\Omega$.

With the help of (\ref{n0544}) and (\ref{n0544n}), we can derive more
uniform bounds on $(\rho,\mathbf{v})$. Using (\ref{0404}) and (\ref{n0544n}),
thanks to the norm equivalence on $\mathbf{X}_n$ (see (\ref{0506})), we find that
\begin{equation}\label{n0547}\begin{aligned}
G_1(\underline{\rho},{\mathcal{E}}_{\delta,0},n,T)\leq
\rho\leq G_2(\bar{\rho},{\mathcal{E}}_{\delta,0},n,T,\Omega),
\end{aligned}\end{equation}
from which, (\ref{n0544}) and (\ref{0506}), it follows that
\begin{equation}\label{n0548}\begin{aligned}\|\mathbf{v}\|_{C^0(\bar{I}_{n},\mathbf{X}_n)}
\leq G(\underline{\rho},{\mathcal{E}}_{\delta,0},n,T,\Omega)
\end{aligned}\end{equation}
and
\begin{equation}\label{n0549}\begin{aligned}
\|\mathscr{P}(\rho\mathbf{v})\|_{C^0(\bar{I}_{n}, \mathbf{X}_n)}\leq
c(n)\|\rho\mathbf{v}\|_{C^0(\bar{I}_{n},\mathbf{L}^2(\Omega))}\leq
G(\bar{\rho},{\mathcal{E}}_{\delta,0},n,T,\Omega).\end{aligned}\end{equation}

Applying  (\ref{n0548}) to (\ref{0405}), one gets
\begin{equation}\label{nn0548}\begin{aligned}
\| \rho\|_{C^0(\bar{I}_n,H^1(\Omega))}\leq
G(\underline{\rho},\|\rho_0\|_{H^1(\Omega)},{\mathcal{E}}_{\delta,0}, n, T,\Omega).
\end{aligned}\end{equation}
Utilizing (\ref{0508}), (\ref{0509}), (\ref{n054412}),
(\ref{n0547})--(\ref{n0549}) and (\ref{0506}), we obtain from (\ref{n0413nn}) that
\begin{equation}\label{n0555}
\|\partial_t\mathbf{v}\|_{L^2(I_n,\mathbf{X}_n)}\leq
G(\underline{\rho},\bar{\rho},\|\rho_0\|_{H^1(\Omega)},
\|\mf{d}_0\|_{\mf{H}^2(\Omega)},{\mathcal{E}}_{\delta,0},n,T,\Omega).
\end{equation}
Therefore, we have shown the uniform boundedness of $\|({\rho}^{-1},\rho)\|_{\mf{L}^\infty(Q_{T_n})}$,
$\|\rho\|_{C^0(\bar{I}_n,\mathbf{X}_n)}$, $\|\mathbf{v}\|_{C^0(\bar{I}_n,\mathbf{X}_n)}$,
$\|\mathscr{P}(\rho\mathbf{v})\|_{C^0(\bar{I}_{n}, \mathbf{X}_n)}$ and
$\|\partial_t\mathbf{v}\|_{L^2(I_n,\mathbf{X}_n)}$. It remains to show the uniform boundedness of
$\|\mf{d}\|_{L^\infty(\bar{I}_n,\mf{H}^2(\Omega))}$.

Differentiating (\ref{n302}) with respect to $t$,
multiplying the resulting equations by $\partial_t\mf{d}$,
recalling $|\mf{d}|\leq 1$, we integrate by parts to infer that
\begin{eqnarray}
&& \frac{d}{dt}\int_\Omega|\partial_t\mf{d}|^2\mm{d}\mf{x}  \nonumber \\
&& =2\int_\Omega\partial_t\mf{d}\cdot\left(\theta\Delta\partial_t\mf{d}+\theta f_\varepsilon(|\nabla
\mf{d}|^2)\partial_t\mf{d}+\theta f_\varepsilon'(|\nabla \mf{d}|^2)(\nabla \mf{d}:\nabla
\partial_t\mf{d})\mf{d}-\partial_t\mf{v}\cdot\nabla \mf{d}-\mf{v}\cdot\nabla
\partial_t\mf{d}\right)\mm{d}\mf{x}  \nonumber \\
&& =-2\theta\|\nabla\partial_t\mf{d}\|_{\mf{L}^2(\Omega)}^2 +2\int_\Omega[(\mm{div}\partial_t\mf{v} )
\mf{d}\cdot\partial_t\mf{d}+\partial_t\mf{v}\cdot\nabla \partial_t\mf{d}\cdot\mf{d} -\mf{v}\cdot\nabla
\partial_t\mf{d}\cdot\partial_t\mf{d}]\mm{d}\mf{x} \nonumber \\
&& \quad +2\theta\int_\Omega \partial_t\mf{d}\cdot[f_\varepsilon(|\nabla
\mf{d}|^2)\partial_t\mf{d}+f_\varepsilon'(|\nabla\mf{d}|^2)(\nabla \mf{d}:\nabla
\partial_t\mf{d})\mf{d}]\mm{d}\mf{x} \nonumber \\
&& \leq -\theta\|\nabla
\partial_t\mf{d}\|_{\mathbf{L}^2(\Omega)}^2+C(\theta)\big(\|\partial_t\mf{v}\|_{\mathbf{H}^1(\Omega)}^2
+\|\mf{v}\|_{\mf{L}^\infty(Q_T)}^2\|\partial_t\mf{d}\|_{\mathbf{L}^2(\Omega)}^2
+\|\partial_t\mf{d}\|_{\mathbf{L}^2(\Omega)}^2\big)  \nonumber \\[1mm]
&& \quad +C(\theta)\big(\||\partial_t\mf{d}|^2f_\varepsilon(|\nabla \mf{d}|^2)\|_{{L}^1(\Omega)}
+ \||\partial_t\mf{d}|^2|\nabla \mf{d}|^2 f_\varepsilon'(|\nabla \mf{d}|^2)
\|_{{L}^1(\Omega)}\big),  \label{n228}
\end{eqnarray}
where the last term on the right-hand side of (\ref{n228}) can be bounded as follows.

(1) The two-dimensional case:
 noting that $f_\varepsilon(x)=x$ for the 2D case, we make use of  Lemma \ref{Nirenberg},
and H\"older's and Young's inequalities to see that
\begin{equation}\label{n22291}\begin{aligned}
& C(\theta)\big(\||\partial_t\mf{d}|^2f_\varepsilon(|\nabla
\mf{d}|^2)\|_{{L}^1(\Omega)} + \||\partial_t\mf{d}|^2|\nabla
\mf{d}|^2 f_\varepsilon'(|\nabla \mf{d}|^2)\|_{{L}^1(\Omega)}\big) \\
& \leq C(\theta) \||\partial_t\mf{d}||\nabla \mf{d}|\|^2_{{L}^2(\Omega)}
 \leq C(\theta)\|\partial_t\mf{d}
\|_{\mf{L}^4(\Omega)}^2 \|\nabla \mf{d}\|^2_{\mf{L}^4(\Omega)} \\
& \leq C(\theta,\Omega)(\|\nabla \partial_t\mf{d}\|_{\mf{L}^2(\Omega)}\|\partial_t\mf{d}\|_{\mf{L}^2(\Omega)}
+\|\partial_t\mf{d}\|_{\mf{L^2}(\Omega)}^2) \|\nabla \mf{d}\|^2_{\mf{L}^4(\Omega)} \\
& \leq \frac{\theta}{2}\|\nabla\partial_t\mf{d}\|_{\mathbf{L}^2(\Omega)}^2
+\frac{C(\theta,\Omega)}{\theta}\|\partial_t\mf{d}\|_{\mf{L^2}(\Omega)}^2
(1+\|\nabla \mf{d}\|^2_{\mf{L}^4(\Omega)}).
\end{aligned}\end{equation}

(2) The three-dimensional case: recalling
the definition of $f_\varepsilon(x)$ in \eqref{appfun}, we use \eqref{jww0208jww}, and H\"older's
and Young's inequalities to get
\begin{equation}\label{n2229}\begin{aligned}
& C(\theta)\big(\||\partial_t\mf{d}|^2f_\varepsilon(|\nabla\mf{d}|^2)\|_{{L}^1(\Omega)}
+ \||\partial_t\mf{d}|^2|\nabla \mf{d}|^2 f_\varepsilon'(|\nabla
\mf{d}|^2)\|_{{L}^1(\Omega)}\big)  \\[1mm]
& \leq C(\theta)\|\partial_t\mf{d} \|_{\mf{L}^4(\Omega)}^2
(\|f_\varepsilon(|\nabla \mf{d}|^2)\|_{{L}^2(\Omega)}  +\| |\nabla \mf{d}|^2 f_\varepsilon'(|\nabla \mf{d}|^2)
\|_{{L}^2(\Omega)})\\
& \leq C(\theta,\Omega)\varepsilon^{-1}(\|\nabla \partial_t\mf{d}\|_{\mf{L}^2(\Omega)}^\frac{3}{2}
\|\partial_t\mf{d}\|_{\mf{L}^2(\Omega)}^\frac{1}{2}+\|\partial_t\mf{d}\|_{\mf{L}^2(\Omega)}^2)
\|\nabla \mf{d}\|_{\mf{L}^4(\Omega)}  \\
& \leq \frac{\theta}{2}\|\nabla \partial_t\mf{d}\|_{\mathbf{L}^2(\Omega)}^2
+\frac{C(\theta,\varepsilon^{-4},\Omega)}{\theta}\|\partial_t\mf{d}\|_{\mf{L^2}(\Omega)}^2
(1+\|\nabla \mf{d}\|^4_{\mf{L}^4(\Omega)})\quad\mbox{ for }\;\varepsilon\in (0,1).
\end{aligned}\end{equation}

Inserting \eqref{n22291} and (\ref{n2229}) into (\ref{n228}), we conclude that
\begin{equation*}  \begin{aligned}
&\frac{d}{dt}\|\partial_t\mf{d}\|^2_{\mathbf{L}^2(\Omega)}+\frac{\theta}{2}\|\nabla
\partial_t\mf{d}\|_{\mathbf{L}^2(\Omega)}^2   \\
 &\leq \frac{C(\theta,\varepsilon^{8-4N},\Omega)\|\partial_t\mf{d}\|_{\mathbf{L}^2(\Omega)}^2}{\theta}\left(\|\nabla
\mf{d}\|^4_{\mf{L}^4(\Omega)}+ \|\mf{v}\|_{\mf{L}^\infty(Q_T)}^2+1\right)
+C(\theta)\|\partial_t\mf{v}\|_{\mathbf{H}^1(\Omega)}^2,
\end{aligned}\end{equation*}
which, by applying Gronwall's inequality, gives
\begin{equation*}  \begin{aligned}
\|\partial_t\mf{d}\|^2_{\mathbf{L}^2(\Omega)}\leq&\left(
\|\partial_t\mf{d}(0)\|_{\mathbf{L}^2(\Omega)}^2+C(\theta)\|\partial_t\mf{v}\|_{L^2(I_n,\mathbf{H}^1(\Omega))}^2\right)
e^{\frac{C(\theta,\varepsilon^{8-4N},T,\Omega)}{\theta}(\|\nabla
\mf{d}\|^4_{\mf{L}^4(Q_T)}+\|\mf{v}\|_{\mf{L}^\infty(Q_T)}^2+1)}.
\end{aligned} \end{equation*}

Noting that
$$\|\partial_t\mf{d}(0)\|_{\mathbf{L}^2(\Omega)}=\|\theta(\Delta
\mathbf{d}_0+|\nabla \mathbf{d}_0|^2\mf{d}_0)-\mathbf{v}_0\cdot
\nabla \mathbf{d}_0\|_{\mf{L}^2(\Omega)},$$
we use (\ref{n0548}), (\ref{n0555}) and (\ref{n054412}) to arrive at
\begin{equation}\label{nnn0457}\begin{aligned}
\|\partial_t\mf{d}\|^2_{L^\infty(I_n,\mathbf{L}^2(\Omega))}+\|\nabla
\partial_t\mf{d}\|_{L^2(I_n,\mathbf{L}^2(\Omega))}^2 \leq &
G(\underline{\rho},\bar{\rho},\|\rho_0\|_{H^1(\Omega)},
\|\mf{d}_0\|_{\mf{H}^2(\Omega)},{\mathcal{E}}_{\delta,0},n,\varepsilon^{8-4N},T,\Omega).
\end{aligned}\end{equation}

Recalling $|\mf{d}|\leq 1$, from (\ref{n302}) we get
\begin{equation}\label{njww0348}\begin{aligned}
 \theta^2\int_\Omega|\Delta\mf{d}|^2\mm{d}\mf{x}\leq &
 3\theta^2\int_\Omega|\nabla\mf{d}|^4\mm{d}\mf{x}
+3\int_\Omega|\partial_t\mf{d}|^2\mm{d}\mf{x}+3\int_\Omega|\mf{v}|^2
|\nabla\mf{d}|^2\mm{d}\mf{x}.
\end{aligned}\end{equation}
Similar to the derivation of (\ref{jww0219}), the first term on the right-hand side of (\ref{njww0348})
can be bounded as follows.
\begin{equation}\label{n0459}\begin{aligned}
 & 3\theta^2\int_\Omega|\nabla  \mf{d}|^4\mm{d}\mf{x}\leq 3\theta^2\tilde{C}_2(\Omega)
 \tilde{C}_0\epsilon_0\| \Delta \mf{d}\|_{\mathbf{{L}}^{2}(\Omega)}^{2}
 +C(\Omega)(\|\nabla\mf{d}\|_{\mathbf{{L}}^2(\Omega)}^2+
\|\mf{d}-{\mf{e}}_N\|_{\mathbf{L}^2(\Omega)}^{2}),
\end{aligned}\end{equation}
where the constants $\tilde{C}_2(\Omega)$ and $\tilde{C}_0$ are the same as in \eqref{choosingjww0219}.
Noting that $\tilde{C}_2(\Omega)\tilde{C}_0\epsilon_0\leq ({8\tilde{C}_3})^{-1}\leq 4^{-1}$
by \eqref{choosingjww0219}, using \eqref{123654}, (\ref{n054412}) and (\ref{nnn0457}), we
find from \eqref{njww0348} and (\ref{n0459}) that
\begin{equation}\label{nn0460} \begin{aligned}
 \frac{\theta^2}{4}\int_\Omega|\Delta\mf{d}|^2\mm{d}\mf{x}  \leq& 3\theta^2C(\Omega)
 (\|\nabla\mf{d}\|_{\mathbf{{L}}^2(\Omega)}^2 + 1) \\
& + 3\|\partial_t\mf{d}\|_{\mathbf{L}^2(\Omega)}^2+3\|\mf{v}\|_{\mf{L}^\infty(\Omega)}^2\|\nabla
\mf{d}\|_{\mathbf{L}^2(\Omega)}^2 \\
\leq & G(\underline{\rho},\bar{\rho},\|\rho_0\|_{H^1(\Omega)},
\|\mf{d}_0\|_{\mf{H}^2(\Omega)},{\mathcal{E}}_{\delta,0},n,\varepsilon^{8-4N},T,\Omega).
\end{aligned}\end{equation}
Hence, by virtue of (\ref{jww0218}), (\ref{n054412}), \eqref{nnn0457}, (\ref{nn0460}) and the fact $|\mf{d}|\leq 1$,
\begin{equation}\label{fzu0447}\begin{aligned}
 \|\mf{d}\|_{L^\infty(I_n,\mf{H}^2(\Omega))}\leq
& G(\underline{\rho},\bar{\rho},\|\rho_0\|_{H^1(\Omega)},
\|\mf{d}_0\|_{\mf{H}^2(\Omega)},{\mathcal{E}}_{\delta,0},n,\varepsilon^{8-4N},T,\Omega).
\end{aligned}\end{equation}

The inequalities (\ref{n0547}), (\ref{n0548}), (\ref{nn0548}), (\ref{n0555}) and (\ref{fzu0447})
  furnish the desired estimates which, in combination with Proposition \ref{localtime}, give a
possibility to repeat the above fixed point argument to conclude
that $T_n=T$, and moreover, the global solution $(\rho_n,\mf{v}_n,\mf{d}_n)$
is unique. To end this section, we summarize our previous results on
the global existence and uniqueness of a solution $(\rho_n,\mf{v}_n,\mathbf{d}_n)$ to
the third approximate problem (\ref{n301})--(\ref{n0306}) as follows.
%%%%%%%%%%%%%%%%%%%%%%%%%%%%%%%%%%%%%%%%%%%%%%%%%%%%%%%%%%%%%%%%%%%%%
\begin{pro}\label{thm:0301}
Let the constant $\epsilon_0>0$ (depending on $\Omega$) satisfy \eqref{choosingjww0219},
\begin{equation*}\label{condition0452}
\delta>0,\ \beta>0, \ \varepsilon>0,\mbox{ and }0<\underline{\rho}\leq \bar{\rho}<\infty. \end{equation*}
Assume that $\Omega\subset \mathbb{R}^N$ is a bounded $C^{2,\alpha}$-domain ($\alpha\in(0,1)$),
and the initial data $(\rho_0,\mf{m}_0,\mf{d}_0)$ satisfies
\begin{eqnarray}
&&\label{n0467}1-d_{0N}<\epsilon_0,\quad |\mf{d}_0|=1,\quad \mathbf{d}_0\in\mf{H}^{2}({\Omega}),\\
&&\label{0920n} 0<\underline{\rho}\leq \rho_0\leq\bar{\rho},\quad \rho_0\in W^{1,\infty}(\Omega),
\quad  \mf{v}_0\in \mf{X}_n.
\end{eqnarray}
 Then the third approximate problem (\ref{n301})--(\ref{n0306}) possesses a unique triple
$(\rho_n,\mf{v}_n,\mf{d}_n)$ with the following properties:
\begin{enumerate}
\item[\quad \quad(1)] Regularity:
$\rho_n$ satisfies the same regularity as in
Proposition \ref{pro:0401}, $\mf{v}_n\in C^0(\bar{I},\mf{X}_n)$,
$\partial_t\mf{v}_n\in L^2(I,\mathbf{X}_n)$, $\mf{d}_n$ satisfies the same regularity as in
Proposition \ref{pro:0402} with $T$ in place of $T_d^K$.
%%%%%%%%%%%%%%%%%%%%%
  \item[\quad \quad(2)] $(\rho_n,\mf{v}_n,\mf{d}_n)
  $ solves (\ref{n301}) and (\ref{n302}) a.e. in $Q_T$, and
  satisfies (\ref{n303}) and $(\rho_n,\mf{v}_n,\mf{d}_n)|_{t=0}=(\rho_0,\mf{v}_0,\mf{d}_0)$.
%%%%%%%%%%%%%%%%%%%%%%%%%%%%
\item[\quad \quad(3)] Finite and bounded energy inequalities hold in the 2D case:
\begin{equation}\begin{aligned}\label{fintejjw0457enegy}
\frac{d}{dt}{\mathcal{E}}^n_\delta(t)+\mathcal{F}^n(t)+\int_\Omega
\varepsilon\delta\beta \rho_n^{\beta-2}|\nabla
\rho_n|^2(t)\mm{d}\mf{x}\leq 0\;\;\mbox{ in }\mathcal{D}'(I),
\end{aligned}\end{equation}
and
\begin{equation}\begin{aligned}\label{0458enegy}
{\mathcal{E}}_\delta^n(t)+\int_0^t\left(\mathcal{F}^n(s)+\int_\Omega
\varepsilon\delta\beta \rho_n^{\beta-2}|\nabla
\rho_n|^2(s)\mm{d}\mf{x}\right)\mm{d}s\leq
{\mathcal{E}_{\delta}}(\rho_0,\mf{m}_0,\mf{d}_0):= {\mathcal{E}}_{\delta,0}
\end{aligned}\end{equation}
a.e. in $I$, where ${\mathcal{F}}^n(t):={\mathcal{F}}(\rho_n,\mf{v}_n,\mf{d}_n)$ and
${\mathcal{E}}_\delta^n(t):={\mathcal{E}}_\delta(\rho_n,\mf{m}_n,\mf{d}_n)$
with $\mf{m}_n=\rho_n\mf{v}_n$.
%%%%%%%%%%%%%%%%%%%%%%%%%%%%
  \item[\quad \quad(4)] Additional uniform estimates:
\begin{eqnarray} && \label{n047712}
|\mf{d}|\leq 1 \mbox{ in }\bar{Q}_T,\mbox{ in particular, }\, |\mf{d}|=1\; \mbox{ if }\;|\mf{d}_0|=1
\;\mbox{ for }N=2,  \\[1mm]
&&\label{12134n0478}\|\mf{d}_n-\mf{e}_N\|_{\mf{L}^\infty(\Omega)}< C_0\sqrt{\epsilon_0},
\\&&\label{n04771231}\sup_{t\in I}\|(\mf{d}(t)-\mf{e}_N,\nabla \mathbf{d}(t)
)\|_{\mf{L}^2(\Omega)}+\|\nabla^2 \mf{d}\|_{\mf{L}^2(Q_T)}+\|\nabla
\mf{d}\|_{\mf{L}^4(Q_T)}\leq G({\mathcal{E}}_{\delta,0},\Omega),\\
&& \label{n04771231jww}\|\partial_t\mathbf{d}\|_{L^{4/3}(I,\mathbf{L}^2(\Omega))}
+\|\partial_t\mathbf{d}\|_{L^{2}(I,(\mf{H}^{1}(\Omega))^*)}\leq G( {\mathcal{E}}_{\delta,0},\Omega),\\
&& \label{n04771231jww12}\sup\limits_{t\in I}(\|\sqrt{\rho}{\bf
v}(t)\|_{\mathbf{L}^2(\Omega)}+\|\rho(t)\|_{{L}^\gamma(\Omega)}) +\|\nabla
{\mathbf{v}}\|_{L^2(I,\mf{L}^2(\Omega))}\leq G( {\mathcal{E}}_{\delta,0},\Omega),\\
&&\label{n0477} \sqrt{\varepsilon}\|\nabla \rho_n\|_{\mf{L}^2(Q_T)}
\leq G({\mathcal{E}}_{\delta,0},\delta,\Omega),\\
&& \label{n0478} \|\rho_n\|_{L{^\frac{4\beta}{3}}(Q_T)}
 \leq G({\mathcal{E}}_{\delta,0},\varepsilon,\delta,\Omega),
\end{eqnarray}
(see \cite[ Section 7.7.5.2]{NASII04} for the proof of (\ref{n0477}) and (\ref{n0478})),
where $G$ is a positive constant which is independent of $n$ and nondecreasing in
its arguments, and may depend on $T$. Moreover, if $\varepsilon$ is not explicitly written
in the argument of $G$, then $G$ is independent of $\varepsilon$ as well.
%%%%%%%%%%%%%%%%%%%%%%%%%%%%%%%
\item[\quad \quad(5)]
In particular, if $\Omega=B_R$ with $R\geq 1$, then  $\epsilon_0$ can be chosen to be
independent of the domain $\Omega=B_R$ for any $R\geq 1$, and the constant
$G$ in \eqref{n04771231} and \eqref{n04771231jww12} can be replaced by a constant
$C({\mathcal{E}}_{\delta,0},T,\|\mf{d}_0-\mf{e}_N\|_{\mf{L}^2(\Omega)})$ independent of $B_R$.
\end{enumerate}
\end{pro}

\section{Proof of Theorem \ref{thm:0101}}\label{sec:03}

Once we have established Proposition \ref{thm:0301}, we can obtain
Theorem \ref{thm:0101} by employing the standard three-level approximation
scheme and the method of weak convergence in a manner similar to that in \cite{LPLMTFM98,FENAPHOFJ35801}
for the compressible Naiver-Stokes equations. These arguments have also been
successfully used to establish the existence of weak solutions to other models from fluid dynamics,
see the 2D problem of (\ref{0101})--(\ref{0103}) in \cite{JFJSWDH},
and the 3D Ginzburg-Landau approximation model to (\ref{0101})--(\ref{0103}) in \cite{WDHYCGA} for example.

First we can construct a solution sequence $(\rho_n,\mf{v}_n,\mf{d}_n)$ by Proposition \ref{thm:0301}, using
 the related uniform estimates in Proposition \ref{thm:0301} and standard compactness arguments, we can
   obtain the weak limit
$(\rho_\varepsilon,\mf{v}_\varepsilon,\mf{d}_\varepsilon)$ of the solution sequence
$(\rho_n,\mf{v}_n,\mf{d}_n)$ as $n\rightarrow \infty$, taking subsequences if necessary,
which is a weak solution of the following second approximate problem:
\begin{eqnarray*}
&& %\label{n0507f}
\partial_t\rho_\varepsilon+\mm{div}(\rho_\varepsilon\mf{v}_\varepsilon )
-\varepsilon\Delta \rho_\varepsilon =0\mbox{ in }\mathcal{D}'(Q_T),\\[1mm]
&&%\label{n0508f}
\begin{aligned}
&\partial_t(\rho_\varepsilon\mf{v}_\varepsilon)+
\partial_j(\rho_\varepsilon\mf{v}_\varepsilon v_\varepsilon^j)-\mu
\Delta \mf{v}_\varepsilon-(\mu+\lambda)\nabla \mm{div}\mf{v}_\varepsilon+\nabla A\rho^\gamma_\varepsilon
+\delta\nabla \rho_\varepsilon^\beta \\
&\quad +\nu\mathrm{div}\left(\nabla
\mathbf{d}_\varepsilon\otimes\nabla\mathbf{d}_\varepsilon-\frac{|\nabla
\mf{d}_\varepsilon|^2\mathbb{I}}{2}\right)+ \varepsilon
\nabla \rho_\varepsilon\cdot \nabla \mf{v}_\varepsilon=0\;\mbox{ in }
(\mathcal{D}'(Q_T))^N,\end{aligned}\\[1mm]
&& %\label{n0509f}
\partial_t\mathbf{d}_\varepsilon+\mathbf{v}_\varepsilon\cdot \nabla
\mathbf{d}_\varepsilon= \theta(\Delta \mathbf{d}_\varepsilon+f_\varepsilon(|\nabla
\mathbf{d}_\varepsilon|^2)\mf{d}_\varepsilon),\mbox{ a.e. in }Q_T.
\end{eqnarray*}
 with boundary conditions
\begin{eqnarray}\label{njww0366}\nabla\rho_\varepsilon\cdot\mf{n}|_{\partial\Omega}=0,\quad
\mathbf{v}_\varepsilon|_{\partial \Omega}=\mathbf{0},\quad
(\mf{n}\cdot\nabla\mathbf{d}_\varepsilon)|_{\partial\Omega}=\mathbf{0},\end{eqnarray}
and modified initial data (\ref{n0467})--(\ref{0920n}),
where $\delta>0$, $\beta\geq\max\{\gamma,8\}$, and $\varepsilon>0$. Moreover, the solution
$(\rho_\varepsilon,\mf{v}_\varepsilon,\mf{d}_\varepsilon)$ enjoys the finite and bounded energy inequalities \eqref{fintejjw0457enegy}--(\ref{0458enegy}), and uniform estimates
(\ref{n047712})--(\ref{n0478}).

We proceed to utilize the related uniform estimates and standard compactness arguments to
obtain the weak limit $(\rho_\delta,\mf{v}_\delta,\mf{d}_\delta)$ of the weak solution sequence
$(\rho_\varepsilon,\mf{v}_\varepsilon,\mf{d}_\varepsilon)$ to the second approximate problem
as $\varepsilon\to 0$, taking subsequences if necessary, which is a weak solution of the
following first approximate problem:
\begin{eqnarray*}
&& \label{n0507f}\partial_t\rho_\delta+\mm{div}(\rho_\delta\mf{v}_\delta )
=0\quad\mbox{ in }\mathcal{D}'(Q_T),\\[1mm]
&&\label{n0508f}\begin{aligned}
&\partial_t(\rho_\delta\mf{v}_\delta)+
\partial_j(\rho_\delta\mf{v}_\delta v_\delta^j)-\mu
\Delta \mf{v}_\delta-(\mu+\lambda)\nabla \mm{div}\mf{v}_\delta+\nabla A\rho^\gamma_\varepsilon
+\delta\nabla \rho_\delta^\beta \\
&\quad +\nu\mathrm{div}\left(\nabla
\mathbf{d}_\delta\otimes\nabla\mathbf{d}_\delta-\frac{|\nabla
\mf{d}_\delta|^2\mathbb{I}}{2}\right)=0\quad\mbox{ in }
(\mathcal{D}'(Q_T))^2,\end{aligned}\\[1mm]
&& \label{n0509f}
\partial_t\mathbf{d}_\delta+\mathbf{v}_\delta\cdot \nabla
\mathbf{d}_\delta= \theta(\Delta \mathbf{d}_\delta+|\nabla
\mathbf{d}_\delta|^2\mf{d}_\delta),\quad\mbox{ a.e. in }Q_T
\end{eqnarray*}
with boundary conditions (\ref{njww0366}) and modified initial data
(\ref{n0467})--(\ref{0920n}). Moreover, the solution $(\rho_\delta,\mf{v}_\delta,\mf{d}_\delta)$
enjoys the estimates \eqref{n047712}--\eqref{n04771231jww12} and inequalities \eqref{fintejjw0457enegy}--(\ref{0458enegy}) with $\varepsilon=0$.
Here we remark that it is easy to verify the convergence of $f_\varepsilon(|\nabla\mf{d}_\varepsilon|^2)$
to $|\nabla \mf{d}_\delta|^2$ as $\varepsilon\rightarrow 0$ in three dimensions.

Using the uniform bounds given in (\ref{12134n0478})--(\ref{n04771231jww}) with $\mf{d}_\varepsilon$
in place of $\mf{d}_n$, applying the Arzel\`a-Ascoli theorem and Aubin-Lions lemma,
and taking subsequences if necessary, we deduce that
\begin{eqnarray}\label{strong1}
\mathbf{d}_\varepsilon\rightarrow \mathbf{d}_\delta\;\mbox{ strongly in }\;
C^0(I,\mathbf{L}^2(\Omega))\cap L^p(I, \mf{H}^{1}(\Omega))\cap L^r(I, \mf{W}^{1,q}(\Omega))
\end{eqnarray}
for any $p\in [1,\infty)$, $q\in [1,6)$ and $r\in [1,2)$, which, recalling the definition of $f_\varepsilon$,
implies that (taking subsequences if necessary)
\begin{eqnarray}\label{302}
f_\varepsilon(|\nabla \mathbf{d}_\varepsilon|^2)\mf{d}_\varepsilon\to |\nabla\mathbf{d}_\delta|^2\mf{d}_\delta
\;\;\mbox{ as }\varepsilon\to 0\quad\mbox{ for a.e. }\mf{x}\in \Omega.
\end{eqnarray}
Thus, using Vitali's convergence theorem, and recalling the uniform in $\varepsilon$ boundedness
of $\|\nabla \mf{d}_\varepsilon\|_{\mf{L}^2(Q_T)}$, we infer that
\begin{eqnarray}\label{strong2}\begin{aligned}
f_\varepsilon(|\nabla \mathbf{d}_\varepsilon|^2)\mf{d}_\varepsilon\to |\nabla
\mathbf{d}_\delta|^2\mf{d}_\delta & \;\;\mbox{ strongly in } \mf{L}^r(Q_T)\;\;\mbox{ for any }r\in [1,2),\\
& \mbox{ and weakly in }\mf{L}^2(Q_T).
\end{aligned}\end{eqnarray}
In addition, we have the regularity $\mf{d}_\delta\in L^2(I,{\mf{H}^2(\Omega)})$
and $\partial_t\mathbf{d}_\delta\in {L^{2}(I,(\mf{H}^{1}(\Omega))^*)}$.
In view of \cite[Proposition 7.31]{NASII04}, we get consequently
\begin{eqnarray}\label{n0504regularity}
\mf{d}_\delta\in C^0(I,{\mf{H}^1(\Omega)}).\end{eqnarray}

Finally, we can also  obtain a weak solution $(\rho,\mf{v},\mf{d})$
of the original problem (\ref{0101})--(\ref{0105})  with boundary
conditions ``$
\mathbf{v}|_{\partial \Omega}=\mathbf{0}$ and $
(\mf{n}\cdot\nabla\mathbf{d})|_{\partial\Omega}=\mathbf{0}$",  and modified initial data
(\ref{n0467})--(\ref{0920n}), which is the weak limit as $\delta\to 0$
of the weak solution sequence $(\rho_\delta,\mf{v}_\delta,\mf{d}_\delta)$ of the second approximate problem.
It should be noted that the modified initial energy  in \eqref{0458enegy}--\eqref{n04771231jww12}
can be further chosen to be independent of $\delta$, in other words,
the term ${\mathcal{E}}_{\delta,0}$ in \eqref{0458enegy}--\eqref{n04771231jww12}
can be replace by a positive constant $\bar{\mathcal{E}_0}:=\sup_{0\leq \delta\leq 1}\{\mathcal{E}_{\delta,0}\}$.
Hence, the weak solution $(\rho,\mf{v},\mf{d})$ enjoys the same estimates as in Theorem \ref{thm:0101}.
Applying an approximation argument to the initial data, the modified initial data (\ref{n0467}) and (\ref{0920n})
can be relaxed to \eqref{jfw0110}--\eqref{n0109}.
Consequently, we can obtain the desired Theorem \ref{thm:0101}.
We refer to \cite{JFJSWDH,WDHYCGA} or \cite{FENAPHOFJ35801,NASII04} for the omitted
details of the proof of the limit process and the renormalized solutions \eqref{0109}.

\newcommand\ack{\section*{Acknowledgement}}
\newcommand\acks{\section*{Acknowledgements}}
\acks
  The research of Fei Jiang was supported by NSFC (Grant No. 11101044 and 11271051), the research of Song Jiang
  by the National Basic Research Program under the Grant 2011CB309705 and NSFC (Grant No. 11229101),
  while the research of Dehua Wang by the National Science Foundation under Grant DMS-0906160
 and the Office of Naval Research under Grant N00014-07-1-0668.

% The authors would like to thank the anonymousx referee for invaluable
%suggestions.
%
\renewcommand\refname{References}
\renewenvironment{thebibliography}[1]{%
\section*{\refname}
\list{{\arabic{enumi}}}{\def\makelabel##1{\hss{##1}}\topsep=0mm
\parsep=0mm
\partopsep=0mm\itemsep=0mm
\labelsep=1ex\itemindent=0mm
\settowidth\labelwidth{\small[#1]}%
\leftmargin\labelwidth \advance\leftmargin\labelsep
\advance\leftmargin -\itemindent
\usecounter{enumi}}\small
\def\newblock{\ }
\sloppy\clubpenalty4000\widowpenalty4000
\sfcode`\.=1000\relax}{\endlist}
\bibliographystyle{model1b-num-names}

\begin{thebibliography}{37}
\expandafter\ifx\csname
natexlab\endcsname\relax\def\natexlab#1{#1}\fi
\providecommand{\bibinfo}[2]{#2} \ifx\xfnm\relax
\def\xfnm[#1]{\unskip,\space#1}\fi
%Type = Book
%\bibitem[{Adams and John(2005)}]{ARAJJFF}
%\bibinfo{author}{R. A. Adams}, \bibinfo{author}{J.~John},
%  \bibinfo{title}{{Sobolev Space}}, \bibinfo{publisher}{Academic Press: New
%  York}, \bibinfo{year}{2005}.

 %Type = Article
%\bibitem[{Ding and Wen(2012)}]{DSJHJRXFG}
%\bibinfo{author}{S.~Ding}, \bibinfo{author}{J.R. Huang}, \bibinfo{author}{F. G. Xia},
%\bibinfo{title}{{A remark on global existence of strong solution
%for incompressible hydrodynamic flow of liquid crystals with vacuum}},
%  \bibinfo{journal}{Preprint}  \bibinfo{year}{2012}.

\bibitem[{Ding and Wen(2013)}]{DSJHJRHYWRZZ}
\bibinfo{author}{S.~Ding}, \bibinfo{author}{J.R. Huang}, \bibinfo{author}{F. G. Xia},\bibinfo{author}{ H. Y. Wen},\bibinfo{author}{R. Z. Zi},
\bibinfo{title}{{Incompressible limit of the compressible nematic liquid
crystal flow}},
  \bibinfo{journal}{J. Funct. Anal.   http://dx.doi.org/10.1016/j.jfa.2013.01.011}  \bibinfo{year}{(2013)}.

 %Type = Article
\bibitem[{Ding et~al.(2011)Ding, Lin, Wang and Wen}]{DSLJWCWH}
\bibinfo{author}{S.~Ding}, \bibinfo{author}{J.~Lin}, \bibinfo{author}{C.~Wang},
  \bibinfo{author}{H.~Wen}, \bibinfo{title}{{Compressible hydrodynamic flow of
  liquid crystals in 1-D}}, \bibinfo{journal}{Discrete and Continuous Dynamical
  Systems--Series B} \bibinfo{volume}{15} (\bibinfo{year}{2011}),
  \bibinfo{pages}{357--371}.
%Type = Article
\bibitem[{Ding and Wen(2011)}]{DSJWHYSN}
\bibinfo{author}{S.~Ding}, \bibinfo{author}{H.~Wen}, \bibinfo{title}{{Solutions
  of incompressible hydrodynamic flow of liquid crystals}},
  \bibinfo{journal}{Nonlinear Analysis: Real World Applications}
  \bibinfo{volume}{12} (\bibinfo{year}{2011}), \bibinfo{pages}{1510--1531}.

\bibitem[{Evans(1998)}]{ELGP}
\bibinfo{author}{L.C. Evans}, \bibinfo{title}{{Partial Differential Equations}},
  \bibinfo{publisher}{American Mathematical Society},
  \bibinfo{address}{Providence, RI}, \bibinfo{year}{1998}.

\bibitem[{Feireisl et~al.(2001)Feireisl, Novotn{\`y} and
  Petzeltov{\'a}}]{FENAPHOFJ35801}
\bibinfo{author}{E.~Feireisl}, \bibinfo{author}{A.~Novotn{\`y}},
  \bibinfo{author}{H.~Petzeltov{\'a}}, \bibinfo{title}{{On the existence of
  globally defined weak solutions to the Navier-Stokes equations}},
  \bibinfo{journal}{Journal of Mathematical Fluid Mechanics}
  \bibinfo{volume}{3} (\bibinfo{year}{2001}), \bibinfo{pages}{358--392}.

  \bibitem{HW1} J. Hineman, C. Wang,
  Well-posedness of Nematic liquid crystal flow in $L^3_{{uloc}}(R^3)$,
    arXiv:1208.5965 [math.AP].

%Type = Article
\bibitem[{Hu and Wu(2012)}]{HXWHGS}
\bibinfo{author}{X.~Hu}, \bibinfo{author}{H.~Wu}, \bibinfo{title}{{Global
  solution to the three-dimensional compressible flow of liquid crystals}},
  \bibinfo{journal}{arXiv:1206.2850v1 [math.AP]}
  (\bibinfo{year}{2012}).
%Type = Article
\bibitem[{Huang et~al.(2012{\natexlab{a}})Huang, Y. and Wen}]{HTWCYWHY}
\bibinfo{author}{T.~Huang}, \bibinfo{author}{C. Wang}, \bibinfo{author}{H.
  Wen}, \bibinfo{title}{{Blow up criterion for compressible nematic liquid
  crystal flows in dimension three}}, \bibinfo{journal}{Arch. Rational Mech.
  Anal.} \bibinfo{volume}{204} (\bibinfo{year}{2012}{\natexlab{a}}),
  \bibinfo{pages}{285--311}.
%Type = Article
\bibitem[{Huang et~al.(2012{\natexlab{b}})Huang, Y. and Wen}]{HTWCYWWHYSS}
\bibinfo{author}{T.~Huang}, \bibinfo{author}{C. Wang}, \bibinfo{author}{H. Wen},
  \bibinfo{title}{{Strong solutions of the compressible nematic liquid
  crystal flow}}, \bibinfo{journal}{Journal of Differential Equations}
  \bibinfo{volume}{252} (\bibinfo{year}{2012}{\natexlab{b}}),
  \bibinfo{pages}{2222--2265}.

 %Type = Article
\bibitem[{Jiang and Jiang(2009)}]{JFJSWDH}
\bibinfo{author}{F.~Jiang}, \bibinfo{author}{S.~Jiang}, \bibinfo{author}{D. H.~Wang},
\bibinfo{title}{{ Global weak solutions to the equations of compressible flow of
nematic liquid crystals in two dimensions}},
  \bibinfo{journal}{arXiv:1210.3565v1 [math.AP] 12 Oct 2012}.

%Type = Article
\bibitem[{Jiang and Tan(2009)}]{JFTZOGWS}
\bibinfo{author}{F.~Jiang}, \bibinfo{author}{Z.~Tan},
\bibinfo{title}{{Global weak solution to the flow of liquid crystals system}},
  \bibinfo{journal}{Math. Meth. Appl. Sci.} \bibinfo{volume}{32}
  (\bibinfo{year}{2009}), \bibinfo{pages}{2243--2266}.

%Type = Article
\bibitem[{Jiang and Zhang(2001)}]{JSZPO}
\bibinfo{author}{S.~Jiang}, \bibinfo{author}{P.~Zhang}, \bibinfo{title}{{On
  spherically symmetric solutions of the compressible isentropic Navier-Stokes
  equations}}, \bibinfo{journal}{Commun. Math. Phys.} \bibinfo{volume}{215}
  (\bibinfo{year}{2001}), \bibinfo{pages}{559--581}.
%Type = Article
\bibitem[{Lei et~al.(2012)Lei, Li and Zhang}]{LZLDZXY}
\bibinfo{author}{Z. Lei}, \bibinfo{author}{D.~Li}, \bibinfo{author}{X.
  Zhang}, \bibinfo{title}{A new proof of global wellposedness of liquid
  crystals and heat harmonic maps in two dimensions},
  \bibinfo{journal}{arXiv:1205.1269v2 [math.AP]}
  (\bibinfo{year}{2012}).
%Type = Article
\bibitem[{Li and Xu(2012)}]{LJXZHZJWG}
\bibinfo{author}{J.~Li}, \bibinfo{author}{Z. Xu, J. Zhang},
  \bibinfo{title}{{Global well-posedness with large oscillations and vacuum to
  the three-dimensional equations of compressible nematic liquid crystal flows}},
  \bibinfo{journal}{arXiv:1204.4966v1 [math.AP] }
  (\bibinfo{year}{2012}).
%Type = Article

\bibitem{LW-JDE2012}
X. Li, D. Wang, \bibinfo{title}{Global solution to the incompressible flow of liquid crystals},
 \bibinfo{journal}{J. Differential Equations} \bibinfo{volume}{252}
 (\bibinfo{year}{2012}), \bibinfo{pages}{745--767}.

\bibitem[{Li and Wang(2012)}]{LXWDHG}
\bibinfo{author}{X.~Li}, \bibinfo{author}{D.~Wang}, \bibinfo{title}{Global
  strong solution to the density-dependent incompressible flow of liquid
  crystals}, \bibinfo{journal}{to appear in Transactions of AMS. arXiv:1202.1011v1 [math.AP]}
  (\bibinfo{year}{2012}).
%Type = Article
\bibitem[{Lin(1989{\natexlab{a}})}]{LFHNC}
\bibinfo{author}{F.~Lin}, \bibinfo{title}{Nonlinear theory of defects in
  nematic liquid crystal: phase transition and flow phenomena},
  \bibinfo{journal}{Comm. Pure Appl. Math.} \bibinfo{volume}{42}
  (\bibinfo{year}{1989}{\natexlab{a}}), \bibinfo{pages}{789--814}.
%Type = Article
%\bibitem[{Lin(1989{\natexlab{b}})}]{LFHNC4}
%\bibinfo{author}{F.~Lin}, \bibinfo{title}{{Nonlinear theory of defects in
%  nematic liquid crystals: Phase transition and flow phenomena}},
%  \bibinfo{journal}{CPAM} \bibinfo{volume}{42}
%  (\bibinfo{year}{1989}{\natexlab{b}}), \bibinfo{pages}{789--814}.
%Type = Article
\bibitem[{Lin et~al.(2011)Lin, Lin and Wang}]{LFHLJYWCY}
\bibinfo{author}{F.~Lin}, \bibinfo{author}{J.~Lin}, \bibinfo{author}{C.~Wang},
  \bibinfo{title}{Liquid crystal flows in two dimensions},
  \bibinfo{journal}{Arch. Rational Mech. Anal} \bibinfo{volume}{197}
  (\bibinfo{year}{2011}), \bibinfo{pages}{297--336}.
%Type = Article
\bibitem[{Lin and Liu(1995{\natexlab{a}})}]{LFHLCNC5}%LFHLCNCX}
\bibinfo{author}{F.~Lin}, \bibinfo{author}{C.~Liu},
  \bibinfo{title}{Nonparabolic dissipative systems modeling the flow of liquid
  crystals}, \bibinfo{journal}{Comm. Pure Appl. Math.} \bibinfo{volume}{XLV
  III} (\bibinfo{year}{1995}{\natexlab{a}}), \bibinfo{pages}{501--537}.
%Type = Article
%\bibitem[{Lin and Liu(1995{\natexlab{b}})}]{LFHLCNC5}
%\bibinfo{author}{F.~Lin}, \bibinfo{author}{C.~Liu},
 % \bibinfo{title}{{Nonparabolic dissipative systems modeling the flowof liquid
%  crystals}}, \bibinfo{journal}{CPAM} \bibinfo{volume}{XLVIII}
%  (\bibinfo{year}{1995}{\natexlab{b}}) \bibinfo{pages}{501--537}.
%Type = Article
\bibitem[{Lin and Liu(1996)}]{LFHLCPD2}
\bibinfo{author}{F.~Lin}, \bibinfo{author}{C.~Liu}, \bibinfo{title}{Partial
  regularities of the nonlinear dissipative systems modeling the flow of liquid
  crystals}, \bibinfo{journal}{Discrete Cont. Dyn. S.} \bibinfo{volume}{2}
  (\bibinfo{year}{1996}), \bibinfo{pages}{1--23}.
%Type = Article
\bibitem[{Lin and Ding(2012)}]{LJYDSJO}
\bibinfo{author}{J.Y. Lin}, \bibinfo{author}{S.J. Ding}, \bibinfo{title}{On the
  well-posedness for the heat flow of harmonic maps and the hydrodynamic flow
  of nematic liquid crystals in critical spaces},
  \bibinfo{journal}{Math. Meth. Appl. Sci.} \bibinfo{volume}{35} (\bibinfo{year}{2012}),
  \bibinfo{pages}{158--173}.
%Type = Book
\bibitem[{Lions(1998)}]{LPLMTFM98}
\bibinfo{author}{P.~Lions}, \bibinfo{title}{{Mathematical Topics in Fluid
  Mechanics: Compressible Models}}, \bibinfo{publisher}{Oxford University
  Press}, \bibinfo{address}{Oxford}, \bibinfo{year}{1998}.
%Type = Article
\bibitem[{Liu and Zhang(2009)}]{LXZZLC3121}
\bibinfo{author}{X.~Liu}, \bibinfo{author}{Z.~Zhang}, \bibinfo{title}{$L^p$
  existence of the flow of liquid crystals system}, \bibinfo{journal}{Chinese
  Ann. Math.} \bibinfo{volume}{30A} (\bibinfo{year}{2009}), \bibinfo{pages}{1--20}.
%Type = Article
\bibitem[{Liu and Qing(2011)}]{LXGQJE}
\bibinfo{author}{X.G. Liu}, \bibinfo{author}{J.~Qing},
  \bibinfo{title}{Existence of globally weak solutions to the flow of
  compressible liquid crystals system}, \bibinfo{journal}{Preprint}
  (\bibinfo{year}{2011}).
%Type = Article
\bibitem[{Liu and Qing(2011)}]{Nirenberg}
\bibinfo{author}{L. Nirenberg},
  \bibinfo{title}{On elliptic partial differential equations}, \bibinfo{journal}{
  Estratto dagli Annali della Scuola Normale Superiore di Pisa Serie III}
\bibinfo{volume}{XIII. Fasc. II} (\bibinfo{year}{1959}).

%Type = Book
\bibitem[{Novotn{\`y} and Stra{\v{s}}kraba(2004)}]{NASII04}
\bibinfo{author}{A.~Novotn{\`y}}, \bibinfo{author}{I.~Stra{\v{s}}kraba},
  \bibinfo{title}{{Introduction to the Mathematical Theory of Compressible
  Flow}}, \bibinfo{publisher}{Oxford University Press, Oxford},
  \bibinfo{year}{2004}.
%Type = Book
\bibitem[{Temam(1984)}]{NSRT1}
\bibinfo{author}{R.~Temam}, \bibinfo{title}{Navier--Stokes Equations: Theory
  and Numerical Analysis}, \bibinfo{publisher}{R. Temam},
  \bibinfo{address}{Amsterdam}, \bibinfo{year}{1984}.
  %Type = Article
\bibitem[{Suen and Hoff(2012)}]{SADHG}
\bibinfo{author}{A. Suen}, \bibinfo{author}{D. Hoff},
  \bibinfo{title}{ Global low-energy weak solutions of the equations
of three-dimensional compressible magnetohydrodynamics},
\bibinfo{journal}{Arch. Rational Mech. Anal.},
 \bibinfo{volume}{205}(\bibinfo{year}{2012}),
 \bibinfo{pages}{27--58}.
%Type = Article
\bibitem[{Wang(2011)}]{WCYWA}
\bibinfo{author}{C. Y. Wang}, \bibinfo{title}{{Well-posedness for the heat flow
  of harmonic maps and the liquid crystal flow with rough initial data}},
  \bibinfo{journal}{Arch. Rational Mech. Anal.} \bibinfo{volume}{200}
  (\bibinfo{year}{2011}),  \bibinfo{pages}{1--19}.
  %Type = Article
%\bibitem[{Wang(2011)}]{WCYWAHF}
%\bibinfo{author}{C.Y. Wang}, \bibinfo{title}{{Heat Flow of Harmonic Maps Whose
%Gradients Belong to $L^n_xL_t^\infty$}},
%  \bibinfo{journal}{Arch. Rational Mech. Anal.} \bibinfo{volume}{188}
%  (\bibinfo{year}{2008}),  \bibinfo{pages}{351--369}.
%Type = Article
\bibitem[{Wang and Yu(2012)}]{WDHYCGA}
\bibinfo{author}{D.~Wang}, \bibinfo{author}{C.~Yu}, \bibinfo{title}{{ Global
  weak solution and large-time behavior for the compressible flow of liquid
  crystals}}, \bibinfo{journal}{Arch. Rational Mech. Anal.}
  \bibinfo{volume}{204} (\bibinfo{year}{2012}),  \bibinfo{pages}{881--915}.
%Type = Article
%\bibitem[{Wen and Ding(2011)}]{WHYDSJSN}
%\bibinfo{author}{H.~Wen}, \bibinfo{author}{S.~Ding}, \bibinfo{title}{{
%  Solutions of incompressible hydrodynamic flow of liquid crystals}},
%  \bibinfo{journal}{Nonlinear Analysis: Real World Applications}
 % \bibinfo{volume}{12} (\bibinfo{year}{2011}) \bibinfo{pages}{1510--1531}.

%Type = Article
\bibitem[{Wu and Tan(2012)}]{WGTZGL}
\bibinfo{author}{G. C. Wu}, \bibinfo{author}{Z. Tan},
\bibinfo{title}{{Global low-Energy weak solution and large-time behavior for the
compressible flow of liquid crystals}}, \bibinfo{journal}{Preprint}
 \bibinfo{year}{2012}.

\end{thebibliography}

\end{document}